\documentclass[11pt]{article}
\usepackage{amsmath,amssymb}

\newcommand{\prf}{\noindent{\bf Proof.}\ }
\newcommand{\Tbar}{\ensuremath{\overline{\mathcal{T}}}}
\newcommand{\Mbar}{\overline{\mathcal{M}}}
\newcommand{\trip}{|||}
\newcommand{\grad}{\mbox{grad\,}}

\newtheorem{definition}{Definition}

\newtheorem{theorem}[definition]{Theorem}
\newtheorem{corollary}[definition]{Corollary}
\newtheorem{proposition}[definition]{Proposition}

\begin{document}

\title{Geometry of the Weil-Petersson completion of Teichm\"{u}ller space}
\author{Scott A. Wolpert\thanks{Research supported in part by NSF Grants
DMS-9800701.}}
\date{August 21, 2003}
\maketitle

\section{Introduction}
Let $R$ be a Riemann surface of genus $g$ with $n$ punctures, $3g-3+n>0$, and
$\mathcal T$ the Teichm\"{u}ller space  of $R$.  The Weil-Petersson (WP)
metric for $\mathcal T$  is a  K\"{a}hler metric with negative sectional
curvature \cite{Ahsome,Royicm, Trcurv, Wlchern}.  With the WP metric $\mathcal
T$ is a unique geodesic space \cite{Wlnielsen}: for each pair of points there
is a unique distance-realizing joining curve.  The {\em augmented
Teichm\"{u}ller space} $\Tbar$, a stratified non locally compact space, is the
space of marked noded Riemann surfaces and is a bordification of $\mathcal T$
in the style of Baily-Borel, \cite{Abdegn, Bersdeg}.  For $(g,n)=(1,1)$,
$\Tbar$ is the bordification $\mathbb H\cup\mathbb Q$ of the upper half-plane
with the horoball-neighborhood topology.  The augmented Teichm\"{u}ller space
is in fact the WP metric completion of the Teichm\"{u}ller space \cite{Msext}.
The strata of $\Tbar$ are lower-dimensional Teichm\"{u}ller spaces; each
stratum with its natural WP metric isometrically embeds into the completion
$\Tbar$.  

Our purpose is to present a view of the current understanding of the geometry of the WP geodesics 
on $\Tbar$.  The behavior of geodesics in-the-large has significant consequences for the action of the mapping 
class group; see \cite{Brkwpvs, DW2, MW, Wlnielsen, Yam} and Section 7 below.  The  behavior of geodesics is also an important consideration for the harmonic map problem, as well as the study of rigidity of homomorphisms of lattices 
in Lie groups to the mapping class group \cite{Cor, DKW, DW2, GrSc, KShar, Yam}.  Furthermore, the behavior of geodesics is 
 a consideration for the {\em rank} of $\mathcal T$ \cite{BF}.  
We begin by mentioning a collection of recent results \cite{ Brkwp, Brkwpvs, BF,DW2, MW, McM, Yam}.  
Recall that for a hyperbolic surface, the length of the unique geodesic in a prescribed free homotopy class provides a function, the {\em geodesic-length}, on $\Tbar$ valued in $[0,\infty)$.  A general fact is that geodesic-length 
functions are strictly convex along WP geodesics \cite{Wlnielsen}.

The work of C. McMullen provides a prelude \cite{McM}.  Recall that a Bers embedding 
$\beta_S:\mathcal T \rightarrow \mathbf T_S^*\mathcal T$ is a biholomorphic map of the Teichm\"{u}ller space 
to a domain in a cotangent space; from the Nehari estimate the image is bounded independent of $S$ in terms 
of the Teichm\"{u}ller and the WP co-metrics.  Observe for $S_0$ fixed, $-\beta_S(S_0)$  is a section of the cotangent 
bundle $\mathbf T^*\mathcal T$, a differential $1$-form $\theta_{WP}(S)$ on $\mathcal T$.   McMullen 
showed \cite[Thrm. 1.5]{McM} that $d(i\theta_{WP})=\omega_{WP}$ is the WP symplectic form.  An application is a positive lower bound for the WP Rayleigh-Ritz quotient.    He then introduced a smooth modification of the WP metric by including the complex Hessians of the 
small-valued geodesic-length 
functions.  He combined the above and estimates for geodesic-length derivatives to show that the modification is a {\em K\"{a}hler hyperbolic metric} for the moduli space of Riemann surfaces 
that is comparable to the Teichm\"{u}ller metric \cite[Thrm. 1.1]{McM}.  As applications McMullen found: a positive lower bound for the Teichm\"{u}ller Rayleigh-Ritz quotient, a complex submanifold isoperimetric 
inequality, and the alternating sign of the orbifold Euler characteristic for the moduli spaces \cite{McM}.

J. Brock has established important results on the
large-scale behavior of WP distance,\cite{Brkwp}.  Brock considered the metric space, the {\em pants graph}
$C_{\mathbf P}(F)$, having vertices the distinct pants decompositions of $F$
and joining edges of unit-length for pants decompositions differing by a single
elementary move,\cite{Brkwp}.  He showed that the $0$-skeleton of $C_{\mathbf
P}(F)$ is quasi-isometric to  $\mathcal T$ with the WP metric.  In particular by an observation of L. Bers 
there is a constant $L$ such that each hyperbolic surface has a pants decomposition by geodesics of length 
at most $L$.  For a pants decomposition $\mathcal P$, denote by $V(\mathcal P)\subset\mathcal T$ the 
subset of surfaces with the designated decomposition.  The union $\cup_{\mathcal P}V(\mathcal P)$ provides an 
open cover for $\mathcal T$.  Brock found that WP distance records the configuration of the open sets $V(\mathcal P)$ with the $0$-skeleton of  $C_{\mathbf P}(F)$ as the metric model.  An important consequence of Brock's result is the 
correspondence between {\em quasi-geodesics} (quasi length-minimizing paths) on $\mathcal T$ and {\em quasi-geodesics} 
on  $C_{\mathbf P}(F)$.   He further showed for $p,q\in \mathcal T$ that the corresponding quasifuchsian 
hyperbolic three-manifold has convex-core volume comparable to 
$d_{WP}(p,q)$.  At large-scale WP distance and convex-core volume are approximately combinatorially determined.  He also 
showed that the first eigenvalue of the hyperbolic manifold and corresponding Hausdorff 
dimension of the limit set are estimated in terms of WP distance. 

J. Brock and B. Farb used the correspondence to study the 
{\em rank} of $\mathcal T$ in the sense of M. Gromov \cite{BF}.  A notion for the rank of a metric space is the maximal dimension of a quasi-flat, a quasi-isometric embedding of a Euclidean space.  Brock and Farb found that 
$C_{\mathbf P}(F)$ contains quasi-flats of dimension $g-1 + \lfloor \frac{g+n}{2}\rfloor$.  It follows from application of Brock's quasi-isometry that the WP rank is likewise bounded.  
Gromov-hyperbolic metric spaces have rank one and thus the bound  provides for $\dim \mathcal T >2$, 
that $\mathcal T$ is not Gromov-hyperbolic \cite[Thrm. 1.1]{BF}.  The authors further found for $\dim \mathcal T\le 2$ that $C_{\mathbf P}(F)$ and  thus $\mathcal T$ are Gromov-hyperbolic \cite[Thrm. 5.1]{BF}. Yamada and M. Bestvina had also
considered the maximal dimension of a flat, \cite{Yamcom}.  Z. Huang has discovered further new asymptotic flatness \cite{Zh}.  Variation of independent plumbing parameters $t$ prescribes planes with WP curvature $O((-\log |t|)^{-1})$.

W. Ballman and P. Eberlein posed a group-theoretic notion of the rank \cite{Ivrk}.  For discrete cofinite isometry groups of complete simply connected Riemannian manifolds with non positive curvature bounded from below the Ballman-Eberlein notion coincides with the geometrically defined rank.   N. Ivanov has shown that mapping class groups have rank one \cite{Ivrk}.  N. Ivanov and independently B. Farb, A. Lubotzky and Y. Minsky further proved that any infinite-order element in the mapping class group has linear growth in the word metric; at least $O(n)$ generators of the group are required to write 
the $n^{th}$ iterate of an element of infinite order \cite{FLM, Ivmcg}. Rank-1 lattices in simple Lie 
groups have the $O(n)$ writing-property, while higher-rank lattices do not have the property.

An important discovery of Sumio Yamada was the non refraction of WP geodesics: a
geodesic on $\Tbar$ at most changes strata at its endpoints; see \cite[Thrm. 2]{Yam}, \cite[Lemma 3.6]{DW2} and 
Propositions \ref{refract} and \ref{corners} below. A second important observation was that the strata of $\Tbar$ are geodesically convex.  Yamada refined the original WP expansion of Masur 
\cite{Msext} to present a third-order remainder expansion of the metric in the $C^1$-category.  A key ingredient
was the use of an improved estimate for degenerating families of hyperbolic
metrics.  The considerations were based on the relatively technical work of
Wolf \cite{Wfinf} and Wolf-Wolpert \cite{WWana}.  Yamada used the expansion to study the behavior of geodesics in neighborhoods of the bordification. He considered the WP Levi-Civita connection and one-dimensional harmonic
maps to investigate the non refraction.  Yamada then used the convexity of geodesic-length functions and the negative WP curvature to find that $\Tbar$ is a $CAT(0)$ space; see \cite{DW2}, the attribution to B. Farb in \cite{MW} and  Theorem \ref{catzero} below.  He further noted that geodesic convexity of strata is an immediate consequence of the convexity of geodesic-length functions \cite[Thrm. 1]{Yam}.
He applied the statements to give consideration of fixed-points and realizing translation lengths for mapping classes. Yamada also presented that irreducible elements of the mapping class group have positive translation length and a unique axis.  The work has served as an inspiration for the work of Daskalopoulos and Wentworth \cite{DW2}, as well as the author.  

The geometry of $CAT(0)$ spaces is developed in Bridson-Haefliger \cite{BH}. A {\em geodesic triangle} is 
prescribed by a triple of points and a triple of joining length-minimizing curves. A characterization of curvature 
for metric spaces is provided in terms of distance-comparisons for geodesic triangles.   In a $CAT(0)$ space the 
distance and angle measurements for a triangle are bounded by the corresponding measurements for a Euclidean 
triangle with the corresponding edge-lengths \cite[Chap. II.1, Prop. 1.7]{BH}.

G. Daskalopoulos and R. Wentworth gave an independent treatment of the WP
expansion, the non refraction, the $CAT(0)$ result and a more extensive
consideration of actions of mapping classes \cite{DW2}.  The authors obtained a
$C^0$-category expansion by applying the cut-and-paste based estimates for
degenerating families of hyperbolic metrics from \cite{Wlhyp}.  
Scaling considerations were used for the energy of a parameterized curve to
establish non refraction.  The authors proved that irreducible
mapping classes have positive translation length and a unique axis \cite[Thrm. 1.1]{DW2}.  
Previously G. Daskalopoulos, L. Katzarkov and R. Wentworth studied the finite energy equivariant
harmonic map problem for the target $\mathcal T$, \cite{DKW}.  In general a condition on an 
isometric action is required for the existence of an energy minimizing equivariant map.  
In the case of a symmetric space target the action should be {\em reductive}.  For 
$\mathcal T$, the authors \cite{DKW, DW2} propose 
{\em sufficiently large} as the counterpart of the reductive hypothesis.  A subgroup of the mapping
class group is sufficiently large provided it contains two irreducible mapping classes acting with distinct
fixed points on the space of projective measured foliations.  Daskalopoulos and Wentworth 
established \cite[Thrm. 6.2]{DW2} divergence of the axes for two as above independent irreducible mapping
classes.  The authors applied their considerations and studied equivariant maps from universal covers of finite volume complete Riemannian manifolds with finitely generated fundamental groups.  They showed that if there is a finite energy
map with sufficiently large image of the fundamental group, then there is a finite energy equivariant harmonic map \cite[Cor. 1.3]{DW2}.

B. Farb and H. Masur established general higher rank superrigidity for the mapping class group as image.  For an irreducible lattice in a semisimple Lie group of $\mathbb R$-rank at least two, a homomorphism to the mapping class group has finite 
image \cite[Thrm. 1.1]{FM}.  The authors also considered homomorphisms from $SL_n(\mathbb Z)$ to the group of
homeomorphism of a surface.  They showed that all homomorphisms are trivial for $n$ greater than an explicit bound in
the genus.

H. Masur and M. Wolf established the WP-analogue of H. Royden's celebrated result: for
$3g-3+n>1$ and $(g,n)\ne (1,2)$, every WP isometry of $\mathcal T$ is induced
by an element of the extended mapping class group.  They considered the asymptotic WP geometry to 
reduce the matter to considering the restriction of an isometry to $\Tbar - \mathcal T$. In particular an isometry
of $\mathcal T$ extends to the completion $\Tbar$; an isometry of $\Tbar$ preserves the strata structure and following
an approach of N. Ivanov agrees with a mapping class on the maximally noded surfaces.  They then established that 
the set of maximally noded surfaces forms a uniqueness set for WP isometries \cite{MW}. 

Brock has also studied the family of WP geodesic rays based at a point, the WP visual sphere, \cite{Brkwpvs}. Rays are considered with the topology of convergence of initial segments.  He established that the action of the mapping class group does not extend continuously to an action on the WP visual spheres, and that the rays to noded surfaces
are dense in the visual spheres.  An additional discovery was that convergence of initial segments in general 
does not provide for convergence of entire rays; see \cite{Brkwpvs} and Section 7 below.

The purpose of this paper is to continue the study in detail of the geometry of WP geodesics on $\Tbar$.  We provide an independent treatment of the WP expansion based on the less technical approach of \cite{Wlhyp}.  We then use 
the opportunity to give a 
range of new applications including: a thorough treatment of the strata structure, a classification of locally
Euclidean subspaces of $\Tbar$, for the Masur-Wolf theorem a new proof based on a convex hull property, and a 
classification of limits of WP geodesics.

We find that $\Tbar$ is  a {\em stratified} unique geodesic space with the strata intrinsically 
characterized by the metric geometry (see Theorem
\ref{onestrata}), \cite{Yam}.  For a reference surface $F$ and $C(F)$, the
partially ordered set {\em the complex of curves},  consider $\Lambda$ the
natural labeling function from $\Tbar$ to $C(F)\cup\{\emptyset\}$.  For a
marked noded Riemann surface $(R,f)$ with $f:F\rightarrow R$, the labeling
$\Lambda((R,f))$ is the simplex of free homotopy classes on $F$ mapped to the
nodes on $R$.  The level sets of $\Lambda$ are the strata of $\Tbar$.   The
unique WP geodesic $\widehat{pq}$ connecting $p,q\in \Tbar$ is contained in
the closure of the stratum with label $\Lambda(p)\cap\Lambda(q)$ (see Theorem
\ref{onestrata}). The open segment $\widehat{pq}-\{p,q\}$ is a solution of the
WP geodesic differential equation on the stratum with label
$\Lambda(p)\cap\Lambda(q)$.   For a point $p$, the stratum with label $\Lambda(p)$ is the 
union of the open geodesic segments containing the point (see Theorem \ref{onestrata}). 

The central consideration is the expansion of the WP metric in a neighborhood of a
point of a positive codimension $m$ stratum $\mathcal S$.  For $s$ a general
multi-index local coordinate for $\mathcal S$ and $t$ a plumbing construction
multi-index parameter for the transverse to $\mathcal S$, we show for the
multi-index parameter $r=(-\log |t|)^{-1/2}$ the following expansion for the
{\em metric symmetric-tensor} (see Corollary \ref{wpnormal}) 
$$
dg^2_{WP}(s,t)\,=\,\bigl(dg^2_{WP}(s,0)\,+\,\pi^3\sum\limits^m_{k=1}(4dr_k^2+r_k^6d\arg^2t_k)\bigr)\,(1+O(\|r\|^3)).
$$
In particular along $\mathcal S$ the WP metric to third-order remainder is a
{\em product-metric} of the WP metric of $\mathcal S$ and metrics
$(4dr^2+r^6d\arg^2t)$, one for each $t$ parameter. The product-structure with
higher-order remainder suggests the isometric embedding of $\mathcal S$ into
$\Tbar$.  In the transverse direction to $\mathcal S$ the metric is {\em
modeled} by the surface of revolution about the $x$-axis of $y=(x/2)^3$.  The
third-order remainder suggests {\em higher-order flatness} for the normal
along $\mathcal S$.  We combine the above expansion, the
rescaling argument for metric spaces and an elementary quadratic inequality to
establish the non refraction of geodesics (see Propositions \ref{refract} and
\ref{corners}).  

Beyond $CAT(0)$, there are important applications for the above expansion. We
are able to combine the {\em flat triangle lemma} of A. D. Alexandrov
\cite{BH} and Theorem \ref{onestrata} to study the locally Euclidean isometric
subspaces (flats) of $\Tbar$.  A classification is established, and it is
found that the maximal dimensional flats are submanifolds of: a product of
Teichm\"{u}ller spaces of  $g$ once-punctured tori and
$\lfloor\frac{g+n}{2}\rfloor -1$ four-punctured spheres (see Proposition
\ref{flatness}). The result is consistent with the conjecture of Brock-Farb
regarding the rank (the maximal dimension of a quasi-isometric embedding of a
Euclidean space) of the WP metric, \cite{BF}.
Following a suggestion of Brock, the considerations also provide that for
$\dim \mathcal T > 2$ the WP metric is not {\em Gromov-hyperbolic}.  Flat
geodesic triangles in $\Tbar-\mathcal T$ are uniformly approximated by
geodesic triangles in $\mathcal T$.

We also investigate applications of the Brock result \cite{Brkwpvs} that the geodesic rays from a point of $\mathcal T$ 
to the noded Riemann surfaces have initial tangents dense in the initial tangent space.  We generalize the result 
and show that the geodesics connecting maximally noded Riemann surfaces have tangents dense in the tangent bundle 
of $\mathcal T$ (see Corollary \ref{density}).  
An immediate consequence is that $\Tbar$ is the closed WP convex hull of the subset of maximally noded Riemann surfaces 
(see Corollary \ref{convexhull}).  The maximally noded Riemann surfaces play a basic role for the WP $CAT(0)$ geometry.  
In Theorem \ref{WPisom} we combine the convex hull property, the intrinsic nature of the strata structure and the classification of simplicial automorphisms of $C_{\mathbf P}(F)$ to study WP isometries.  A new proof of the Masur-Wolf result is provided: for $3g-3+n>1$ and $(g,n)\ne (1,2)$, every WP isometry of $\mathcal T$ is induced
by an element of the extended mapping class group.

The WP metric is mapping class group invariant.   H. Masur found that the
Deligne-Mumford moduli space of stable curves $\Mbar$ is the WP
quotient-metric completion of the moduli space of Riemann surfaces
\cite{Msext}.  We note that the WP metric for $\Mbar$ is not locally uniquely
geodesic near the compactification divisor of noded Riemann surfaces (see
Proposition \ref{notcatzero}).  A complete, convex subset of a $CAT(0)$ space is the base for an 
{\em orthogonal projection}, \cite[Chap. II.2]{BH}.  The closure of a stratum is complete and convex.
We show that the distance to a stratum $\mathcal S$
has an expansion in terms of the defining {\em geodesic-length functions}.
For a positive codimension $m$ stratum $\mathcal S$, defined by the vanishing
of the geodesic-length sum $\ell=\ell_1+\dots+\ell_m$, the distance to the
stratum has the simple expansion $d(\cdot,\mathcal
S)=(2\pi\ell)^{1/2}+O(\ell^2)$ (see Corollary  \ref{distD}). Furthermore the
vector fields  $\{\grad (2\pi\ell_j)^{1/2}\}$ are close to orthonormal near
$\mathcal S$.

Our final application concerns limits of sequences of geodesics.  We consider
the classification problem (see Proposition \ref{geolims}). We might expect the
compactness of $\Mbar$ to be manifested in the sequential compactness of the
space of geodesics.  But Brock already found that convergence of initial segments in general 
does not provide for convergence of entire rays.
In fact for each sequence of bounded length geodesics
there is a subsequence of mapping class group translates that converges
geometrically (sequences of products of Dehn twists are applied to subsegments
of the geodesics) to a polygonal path, a curve piecewise consisting of
geodesics connecting different strata (see Proposition \ref{geolims}).  Polygonal
paths were first considered by Brock in his investigation of the WP visual
sphere and the action of the mapping class group, \cite[esp. Secs. 4,
5]{Brkwpvs}.  We find that the limit polygonal path is unique
length-minimizing amongst paths joining  prescribed strata.  A simple example
of a polygonal path is presented in the opening of Section \ref{geolim}. We
apply the considerations and show  that a mapping class acting on $\Tbar$
either: has a fixed-point,  or positive translation length realized on a
closed convex set, possibly contained in $\Tbar - \mathcal T$ (see Theorem
\ref{axes}). For {\em irreducible} mapping classes, the positive translation
length is realized on a unique geodesic within $\mathcal T$, \cite{DW2,Yam}.   
      
We begin our detailed considerations in the next section with a summary of the notions
associated with lengths of curves in metric spaces, \cite{BH}.  We also review
the local deformation theory of noded Riemann surfaces, as well as the
specification of Fenchel-Nielsen coordinates and the construction of the
augmented Teichm\"{u}ller space.  In the third section we provide the WP
expansion.  We begin considerations with the exact expansion of the hyperbolic
metrics for the model case $zw=t$.  Then we consider in detail families of
noded Riemann surfaces and their hyperbolic metrics.  Beginning with Masur's
description of families of holomorphic $2$-differentials, we give a simple and
self-contained development of the tangent-cotangent coordinate frame pairing
for the local deformation space and the desired WP expansion.  In the fourth
section we develop the length-minimizing properties of the solutions of the WP
geodesic differential equation on $\mathcal T$.  The considerations extend the
earlier treatment \cite{Wlnielsen}.  In the fifth section we develop the
length-minimizing properties of curves on $\Tbar$, including the non
refraction results and the main theorems.  The labeling function $\Lambda$
serves an important role.  WP length-minimizing curves can be analyzed in
terms of their strata-behavior and geodesics within strata.  WP convexity of
the geodesic-length functions also serves an important role. WP geodesics are
confined by the level sets and sublevel sets of geodesic-length functions. Non
refraction is established by a local rescaling of the metric, and an
application of the strict inequality $((a+b)^2+c^2)^{1/2}<(a^2+c^2)^{1/2}+b$
for positive values.  In the sixth section first we examine the circumstance
for the WP distance between corresponding points of a pair of geodesics not
strictly convex.  Then we consider the locally Euclidean isometric
subspaces of $\Tbar$.  We also consider the distance to a stratum.  In the
final section we consider sequences of geodesics and establish the sequential
compactness, as well as a general classification for geodesic limits.  The
results are applied to study the existence of axes for mapping classes.

I would like to thank Jeffrey Brock for conversations.

\section{Preliminaries}
\label{prelim}

We begin with a summary of the notions associated with lengths of curves in a
metric space.  We  closely follow the exposition of Bridson-Haefliger
\cite{BH} and commend their treatment to the reader.  For a metric space
$(M,d)$ the {\em length} of a curve $\gamma:[a,b]\rightarrow M$ is

$$
L(\gamma)=\sup_{a=t_0\le t_1\le \cdots \le t_n=b}\sum\limits_{j=0}^{n-1}d(\gamma(t_j),\gamma(t_{j+1}))
$$
where the supremum is over all possible partitions with no bound on $n$.  A
curve is {\em rectifiable} provided its length is finite.  The basic
properties of length are provided in \cite[Prop. 1.20]{BH}.  Length is lower
semi continuous for a sequence of rectifiable curves converging uniformly to a
rectifiable curve.  A curve $\gamma:[a,b]\rightarrow M$ is {\em parameterized
proportional to arc-length} provided the length of $\gamma$ restricted to
subintervals $[a,t]\subset [a,b]$ is a linear function of $t$, \cite[Defn.
1.21]{BH}. A space $(M,d)$ is a {\em length space} provided the distance
between each pair of points is equal to the infimum of the length of
rectifiable curves joining the points.  It is an observation that the
completion of a length space is again a length space \cite[Exer. 3.6 (3)]{BH}.
A curve $\gamma:[a,b]\rightarrow M$ is {\em length-minimizing} provided for
all $a\le t\le t'\le b$ that
$L(\gamma\vert_{[t,t']})=d(\gamma(t),\gamma(t'))$; we initially reserve the
word {\em geodesic} for curves which are solutions of the geodesic
differential equation on a Riemannian manifold.  A space with every pair of
points having a (unique) length-minimizing joining curve is a {\em (unique)
geodesic space}.  In a metric space a {\em geodesic triangle} is prescribed by
a triple of points and a triple of joining length-minimizing curves.  A
geodesic triangle can be compared to a triangle in a constant-curvature space
with the corresponding sides having equal lengths \cite[Chap. II.1]{BH}.  A
characterization of curvature for metric spaces is provided in terms of
distance-comparisons for comparison triangles \cite[Chap. II.1]{BH}.

Consider $R$ a Riemann surface with complete hyperbolic metric having finite
area. The homeomorphism type of $R$ is given by its genus and number of
punctures.  Relative to a reference topological surface $F$, the surface $R$
is {\em marked} by  an orientation-preserving homeomorphism $f:F \rightarrow
R$.  Marked surfaces $(R,f)$ and $(R',f')$ are {\em equivalent} provided for
$h:R\rightarrow R', h$ a conformal homeomorphism,  $h\circ f$ is homotopic rel
boundary to $f'$.  The set of equivalence classes of the $F$-marked Riemann
surfaces is the Teichm\"{u}ller space $\mathcal{T}$, \cite{ImTan}.  A
neighborhood of the marked surface $(R,f)$ is given by first specifying smooth
Beltrami differentials $\nu_1,\dots,\nu_m$ spanning the Dolbeault group
$H^{0,1}_{\bar{\partial}}(\bar R,\mathcal{E}((\kappa p_1\cdots p_n)^{-1}))$
for $\kappa$ the canonical bundle of the compactification $\bar{R}$ and
$p_1,\dots,p_n$ the point line bundles for the punctures, \cite{Kod}.  For
$s\in \mathbb{C}^m$ set $\nu(s) = \sum_js_j\nu_j$; for $s$ small there is a
Riemann surface $R^{\nu(s)}$ and a diffeomorphism $\zeta^{\nu(s)}:R\rightarrow
R^{\nu(s)}$ satisfying $\bar{\partial}\zeta^{\nu(s)}=\nu(s)
\partial\zeta^{\nu(s)}$.  The parameterization of marked surfaces $s
\rightarrow (R^{\nu(s)}, \zeta^{\nu(s)}\circ f)$ is a holomorphic local
coordinate  for the Teichm\"{u}ller space $\mathcal{T}$.

The {\em mapping class group} $Mod=Homeo^+(F)/Homeo_0(F)$ is the quotient of
the group of orientation-preserving homeomorphisms of $F$ fixing the punctures
by the subgroup of homeomorphisms  isotopic to the identity.  The {\em
extended mapping class group} is the quotient $Mod^*=Homeo(F)/Homeo_0(F)$.  A
mapping class $[h]$ acts on equivalence classes of marked surfaces by taking
$\{(R,f)\}$ to $\{(R,f\circ h^{-1})\}$.  The action of $Mod$ on $\mathcal{T}$
is by biholomorphic maps; the quotient $\mathcal{M}$ is the {\em moduli space}
of Riemann surfaces.  The holomorphic cotangent space of $\mathcal T$ at the
marked surface $(R,f)$ is $Q(R)\backsimeq \check{H}^0(\bar R,
\mathcal{O}(\kappa^2p_1\cdots p_n))$, the space of integrable holomorphic
quadratic differentials.  A  co-metric for the cotangent spaces of
Teichm\"{u}ller space is prescribed by the Petersson Hermitian pairing
$\int_R\varphi\bar{\psi}(dh^2)^{-1}$ for $\varphi,\psi\in Q(R)$ and $dh^2$ the
$R$-hyperbolic metric, \cite{Ahsome}.  The dual metric is the Weil-Petersson
(WP) metric.  The (extended) mapping classes act on $\mathcal{T}$ as WP
isometries; the WP metric projects to $\mathcal M$.  The WP metric is
K\"{a}hler with negative sectional curvature and holomorphic sectional
curvature bounded away from zero, \cite{Royicm, Trcurv, Wlchern}.  Masur
estimated the metric near the compactification divisor $\mathcal{D}$ of the
moduli space, \cite{Msext}.  His preliminary expansion can be used for
after-the-fact insights: the metric is not complete, \cite{Wlnon}; there is an
almost-product structure at infinity, \cite{Yam}; and there are submanifolds
of $\mathcal T$ that approximate Euclidean space (see the present Section
\ref{appl}). The expansion provides that the WP diameter and volume of
$\mathcal{M}$ are finite.  In \cite{Wlhyp} an improved analysis was presented
for the extension of the WP K\"{a}hler form considered {\em in the sense of
currents}. The $\frac{1}{2\pi^2}$ multiple of the WP K\"{a}hler form is the
pushdown of the square of the curvature of the hyperbolic metric considered on
the vertical line bundle for the fibration of the universal curve
$\overline{\mathcal{C}}$ over $\Mbar$.  The multiple of the K\"{a}hler form is
a nonsmooth characteristic class representative of the Mumford class
$\kappa_1$, \cite{Wlhyp}.

The {\em complex of curves} $C(F)$ is defined as follows.  The vertices of
$C(F)$ are (free) homotopy classes of homotopically nontrivial, nonperipheral,
simple closed curves on $F$.  An edge of the complex consists of a pair of
homotopy classes of disjoint simple closed curves.  A $k$-simplex consists of
$k+1$ homotopy classes of mutually disjoint simple closed curves.  A maximal
set of mutually disjoint simple closed curves, a {\em pants decomposition},
has $3g-3+n$ elements.  Brock has described the large-scale WP geometry of
Teichm\"{u}ller space in terms of the {\em pants graph} $C_{\boldsymbol
P}(F)$, a complex whose vertices are the distinct pants decompositions,
\cite{Brkwp}.    The mapping class group $Mod$ acts on curve complexes and in
particular on $C(F)$.

A free homotopy class $\alpha$ of a closed curve on $F$ determines a
geodesic-length function $\ell_{\alpha}$ on $\mathcal T$.  For a marked
surface $(R,f)$, $\ell_{\alpha}$ is the length of the $R$-hyperbolic metric
geodesic homotopic to $f(\alpha)$.  Geodesic-length functions provide
parameters for the Teichm\"{u}ller space.  Suitable collections provide local
coordinates, \cite{ImTan}.  A collection of free homotopy classes
$\{\alpha_1,\dots,\alpha_q\}$ is {\em filling} provided for a set of
representatives with minimal number of self and mutual intersections that
$F-\cup_j\alpha_j$ is a union of topological discs and punctured discs.  A
filling geodesic-length sum $\mathcal L=\sum_j\ell_{\alpha_j}$ is a proper
function on the Teichm\"{u}ller space.  The differential and the WP gradient
of an $\ell_{\alpha}$ are given by the classical Petersson theta-series for
the geodesic.  In \cite{Wlnielsen} we established that
 the WP Hessian of $\ell_{\alpha}$ is positive-definite:
geodesic-length functions are strictly convex along WP geodesics. The
convexity provides a effective way to bound the WP geometry.

The Fenchel-Nielsen coordinates include geodesic-length functions, as well as
lengths of auxiliary segments, \cite{Abbook,ImTan,MasFN, WlFN}.  A pants
decomposition $\mathcal{P}=\{\alpha_1,\dots,\alpha_{3g-3+n}\}$ decomposes the
topological surface $F$ into $2g-2+n$ components (pants), each homeomorphic to
a sphere with a combination of three discs or points removed.  A marked
Riemann surface $(R,f)$ is likewise decomposed into pants by the geodesics
representing $\mathcal P$.  Each component pants, relative to its hyperbolic
metric, has a combination of three geodesic boundaries and cusps.  For each
component pants the shortest geodesic segments connecting boundaries determine
{\em designated points} on each boundary.  For each geodesic in the pants
decomposition of $R$ a parameter $\tau$ is defined as the displacement along
the geodesic between designated points, one  for each side of the geodesic.
For Riemann surfaces close to an initial reference Riemann surface, the
displacement $\tau$ is simply the distance between the designated points; in
general the displacement is the analytic continuation (the lifting) of the
distance measurement.  For $\alpha$ in $\mathcal P$ define the {\em
Fenchel-Nielsen  angle} by $\theta_{\alpha}=2\pi\tau_{\alpha}/\ell_{\alpha}$.
The Fenchel-Nielsen coordinates for Teichm\"{u}ller space for the
decomposition $\mathcal P$ are
$(\ell_{\alpha_1},\theta_{\alpha_1},\dots,\ell_{\alpha_{3g-3+n}},\theta_{\alpha_{3g-3+n}})$.
The coordinates provide a real analytic equivalence of $\mathcal T$ to
$(\mathbb{R}_+\times \mathbb{R})^{3g-3+n}$,\cite{Abbook, ImTan,WlFN}.

A bordification of Teichm\"{u}ller space is introduced by extending the range
of the Fenchel-Nielsen parameters.  The interpretation of {\em length
vanishing} is the key ingredient.  For $\ell_{\alpha}$ equal to zero, the
angle $\theta_{\alpha}$ is not defined and in place of the geodesic for
$\alpha$ there appears a pair of cusps;  $f$ is now a homeomorphism of
$F-\alpha$ to the (marked) hyperbolic surface $R$ (curves parallel to $\alpha$
map to loops encircling the cusps; see the discussion of nodes in the
following Section).  The parameter space for the pair
$(\ell_{\alpha},\theta_{\alpha})$ is the identification space $\mathbb{R}_{\ge
0}\times\mathbb{R}/\{(0,y)\sim(0,y')\}$.  For the pants decomposition
$\mathcal P$ a frontier set $\mathcal{F}_{\mathcal P}$ is added to the
Teichm\"{u}ller space by extending the Fenchel-Nielsen parameter ranges: for
each $\alpha\in\mathcal{P}$, extend the range of $\ell_{\alpha}$ to include
the value $0$, with $\theta_{\alpha}$ not defined for $\ell_{\alpha}=0$.  The
points of $\mathcal{F}_{\mathcal P}$ parameterize ({\em degenerate}) Riemann
surfaces with each $\ell_{\alpha}=0,\alpha\in\mathcal{P},$ specifying a pair
of cusps.  In particular for a simplex $\sigma\subset\mathcal{P}$, the
$\sigma$-null stratum is $\mathcal{S}(\sigma)=\{R\mid \ell_{\alpha}(R)=0
\mbox{ iff }\alpha\in\sigma\}$.  The frontier set $\mathcal{F}_{\mathcal P}$
is the union of the $\sigma$-null strata for the subsimplices of
$\mathcal{P}$.  Neighborhood bases for points of $\mathcal{F}_{\mathcal P}
\subset\mathcal{T}\cup\mathcal{F}_{\mathcal P}$  are specified by the
condition that for each simplex $\sigma\subset\mathcal P$  the projection
$((\ell_{\beta},\theta_{\beta}),\ell_{\alpha}):
\mathcal{T}\cup\mathcal{S}(\sigma)\rightarrow\prod_{\beta\notin\sigma}(\mathbb{R}_+\times\mathbb{R})\times\prod_{\alpha\in\sigma}(\mathbb{R}_{\ge
0})$  is continuous.    For a simplex $\sigma$ contained in pants
decompositions $\mathcal P$ and $\mathcal P'$ the specified neighborhood
systems for $\mathcal T \cup\mathcal{S}(\sigma)$ are equivalent.  The {\em
augmented Teichm\"{u}ller space} $\Tbar=\mathcal{T}\cup_{\sigma\in
C(F)}\mathcal{S}(\sigma)$ is the resulting stratified topological space,
\cite{Abdegn, Bersdeg}.  $\Tbar$ is not locally compact since no point of the
frontier has a relatively compact neighborhood; the  neighborhood bases are
unrestricted in the $\theta_{\alpha}$ parameters for $\alpha$ a $\sigma$-null.
The action of $Mod$ on $\mathcal T$ extends to an action by homeomorphisms on
$\Tbar$ (the action on $\Tbar$ is not properly discontinuous) and the quotient
$\Tbar/Mod$ is (topologically) the compactified moduli space of stable curves
(see the consideration of $\Mbar$ in the next Section), \cite[see Math. Rev.
56 \#679]{Abdegn}.  Masur noted that the WP metric extends to $\Tbar$ and is
complete on $\Mbar$, \cite[Thrm. 2, Cor. 2]{Msext}.  $\Tbar$ is WP complete
since the quotient $\Mbar$ is compact and each point of $\Tbar$ has a
neighborhood with complete closure.

\section{Expansion of the WP metric about the compactification divisor}

Our purpose is to provide a description of {\em local coordinates} for the
local deformation space of a Riemann surface with nodes. We will present  a
modification of the standard coordinates \cite{Bersdeg,Msext} and use the
formulation to present an improved form of Masur's expansion of the WP metric.
The expansion reveals that for the moduli space of stable curves $\Mbar$,
along the compactification divisor $\mathcal{D}$,  the WP metric behaves to
third-order in distance as a product formed with the WP metric of
$\mathcal{D}$.

The description begins with the {\em plumbing variety} 
$\mathcal{V}=\{(z,w,t)\mid zw=t,\ |z|, |w|, |t|<1\}$.  
The defining function $zw-t$ has differential $z\,dw+w\,dz-dt$.  Consequences
are that $\mathcal{V}$ is a smooth variety, $(z,w)$ are global coordinates,
while $(z,t)$ and $(w,t)$ are not.  Consider the projection
$\Pi:\mathcal{V}\rightarrow D$ onto the $t$-unit disc.  $\Pi$ is a submersion,
except at $(z,w)=(0,0)$;  we can consider $\Pi:\mathcal{V}\rightarrow D$ as a
(degenerate) family of open Riemann surfaces.  The $t$-fibre, $t\ne 0$, is the
hyperbola germ $zw=t$ or equivalently the annulus
$\{|t|<|z|<1,\,w=t/z\}=\{|t|<|w|<1,\,z=t/w\}$.  The $0$-fibre is the
intersection of the unit ball with the union of the coordinate axes in
$\mathbb{C}^2$; on removing the origin the union becomes
$\{0<|z|<1\}\cup\{0<|w|<1\}$.  Each fibre of $\mathcal{V}_0=\mathcal{V}-\{0\}
\rightarrow D$ has a complete hyperbolic metric:

\begin{equation}
\begin{split}
\label{hyp}
\mbox{for }t\ne0,&\mbox{ on }\{|t|<|z|<1\} \mbox{ then}\\
&dh^2_t=\Bigl(\frac{\pi}{\log |t|} \csc \frac{\pi\log |z|}{\log |t|}\Bigl|\frac{dz}{z}\Bigr|\Bigr)^2;\\
\mbox{for }t=0,&\mbox{ on }\{0<|z|<1\}\cup \{0<|w|<1\} \mbox{ then}\\
&dh^2_0=\Bigl(\frac{|d\zeta|}{|\zeta|\log|\zeta|}\Bigr)^2\mbox{ for } \zeta=z,\,w.
\end{split}
\end{equation}
The family of hyperbolic metrics $(dh_t^2)$ is a continuous metric, degenerate
only at the origin, for the vertical line bundle of $\mathcal{V}$.  In
particular we have the elementary expansion

\begin{equation}
\begin{split}
\label{hypexp}
dh^2_t=&\, \Bigl(\frac{|d\zeta|}{|\zeta|\log|\zeta|}\Bigr)^2\,\Bigl(\Theta\csc\Theta\Bigr)^2\mbox{\quad for\quad}\Theta=\frac{\pi\log |z|}{\log |t|}\\
=&\, dh_0^2\,\bigl(1+\frac13\Theta^2+\frac{1}{15}\Theta^4+\dots\bigr).
\end{split}
\end{equation}
The parameter $t$ is a boundary point of the annulus $\{|t|<|z|<1\}$.  The
boundary points $t,\,1$ will be included in the data for gluings.  To describe
the variation of annuli with boundary points, we now specify a quasiconformal
map $\zeta$ from the pointed $t$-annulus to the pointed $t'$-annulus
$\zeta(z)=zr^{\beta(r,t')},\, z=re^{i\theta},$ with
$\frac{\partial\beta}{\partial r}$ compactly supported in the annulus.  The
boundary conditions are $\zeta(1)=1$, and by specification
$t|t|^{\beta(|t|,t')}=t'$.  On differentiating in $t'$ and evaluating at
$(|t|,t)$ we find the boundary condition $t\log |t|\,\dot{\beta}(|t|,t)=1$.
More generally the infinitesimal variation of the map is the vector field
$\dot{\zeta}(z)=z\log r\, \dot{\beta}(r,t)$ for $\dot{\zeta}\,, \dot{\beta}$
the first $t$-derivatives.   The map $\zeta$ varies from the identity and has
Beltrami differential 

\begin{equation}
\label{belt}
\bar{\partial}\dot{\zeta}=\frac{z}{2\bar{z}}\,\frac{\partial}{\partial\log r}(\dot{\beta}(r,0)\log r)\,\frac{\overline{dz}}{dz}.
\end{equation}
For sake of later application we evaluate the pairing with a quadratic
differential $z^{\alpha}\bigl(\frac{dz}{z}\bigr)^2$,
\begin{equation}
\label{pair}
\begin{split}
\int_{\{|t|<|z|<1\}}\bar{\partial}\dot{\zeta}\,z^{\alpha}\bigr(\frac{1}{z}\bigl)^2 dE&=\int_{\{|t|<|z|<1\}}\frac{z^{\alpha}}{2z\bar{z}}\frac{\partial}{\partial \log r}(\dot{\beta}\log r)\,dE \\
\mbox{where for }&\alpha=0 \mbox{, then}\\
&=\,\pi \dot{\beta}\log r\Big|^1_{|t|}=\frac{-\pi}{t},\\
\mbox{and other}&\mbox{wise, then}\\
&=0,\\
\end{split}
\end{equation}
for $dE$ the Euclidean area element and where we have applied the boundary
condition for $\dot{\beta}$; the evaluation involves fixing a normalization
for the Serre duality pairing and agrees with \cite[Prop. 7.1]{Msext}.

We review the description of  {\em Riemann surfaces with nodes},
\cite{Bersdeg,Msext,Wlhyp}.  A Riemann surface with nodes $R$ is a connected
complex space, such that every point has a neighborhood isomorphic to either
the unit disc in $\mathbb{C}$, or the germ at the origin in $\mathbb{C}^2$ of
the union of the coordinate axes.  $R$ is {\em stable} provided each component
of $R-\{nodes\}$ has negative Euler characteristic, i.e. has a hyperbolic
metric.  A regular $q$-differential on $R$ is the assignment of a meromorphic
$q$-differential $\Theta_j$ for each component $R_j$ of $R-\{nodes\}$ such
that: i) each $\Theta_*$ has poles only at the punctures of $R_*$ with orders
at most $q$, and ii) if punctures $p,\,p'$ are paired to form a node then
$Res_p\Theta_*=(-1)^q\,Res_{p'}\Theta_*$, \cite{Bersdeg}.

We review the deformation theory of Riemann surfaces with punctures and then
with nodes.  For a Riemann surface $R$ with hyperbolic metric and punctures
there is a natural {\em cusp coordinate} (with unique germ modulo  rotation)
at each puncture: at the puncture $p$, the coordinate $z$ with $z(p)=0$ and
the hyperbolic metric of $R$ given as $\bigl(\frac{|dz|}{|z|\log
|z|}\bigr)^2$, the germ of the hyperbolic metric for the unit-disc.  If the
surface is uniformized by the upper half-plane with $p$ represented by a
width-one cusp at infinity then $z=e^{2\pi i\zeta}$ for $\zeta$ the
uniformization variable.  Now a deformation neighborhood of the marked surface
$R$ is given by specifying smooth Beltrami differentials $\nu_1,\dots,\nu_m$
spanning the Dolbeault group $H^{0,1}_{\bar{\partial}}(\bar
R,\mathcal{E}((\kappa p_1\cdots p_n)^{-1}))$ for $\kappa$ the canonical bundle
of $\bar{R}$ and $p_1,\dots,p_n$ the point line bundles for the punctures,
\cite{Kod}.  For $s\in \mathbb{C}^m$ set $\nu(s) = \sum_ks_k\nu_k$; for $s$
small there is a Riemann surface $R^{\nu(s)}$ and a diffeomorphism
$\zeta:R\rightarrow R^{\nu(s)}$ satisfying $\bar{\partial}\zeta=\nu(s)
\partial\zeta$.  The family of surfaces $\{R^{\nu(s)}\}$ represents a
neighborhood of the marked Riemann surface in its Teichm\"{u}ller space.  We
showed in \cite[Lemma 1.1]{Wlspeclim} that the Beltrami differentials can be
modified a small amount {\em so that} in terms of each {\em cusp coordinate}
the diffeomorphisms $\zeta^{\hat\nu(s)}$ are simply rotations;
$\zeta^{\hat\nu(s)}$is a hyperbolic isometry in a neighborhood of the cusps;
$\zeta^{\hat\nu(s)}$ cannot be complex analytic in $s$, but is real analytic.
We further note that for $s$ small the $s$-derivatives of $\nu(s)$ and
$\hat{\nu}(s)$ are close. We say that $\zeta^{\hat{\nu}(s)}$ {\em preserves
cusp coordinates}.   The parameterization provides a key ingredient for
obtaining simplified estimates of the degeneration of hyperbolic metrics and
an improved expansion for the WP metric.

We review the plumbing construction for $R$ a Riemann surface with a pair of
punctures $p,\,p'$.  The data is $(U,V, F,G,t)$ where: $U$ and $V$ are
disjoint disc coordinate neighborhoods of $p$ and $p'$; $F: U\rightarrow
\mathbb{C},\,F(p)=0$ and $G:V\rightarrow \mathbb{C},\,G(p')=0$, are coordinate
mappings and $t$ is a sufficiently small complex number.  Pick a constant
$0<c<1$ such that $F(U)$ and $G(V)$ contain the disc $\{|\zeta|<c\}$.  For
$c'<c$ let $R^*_{c'}$ be the open surface obtained by removing from $R$ the
discs $\{|F|\le c'\}\subset U$ and $\{|G|\le c'\}\subset V$.  Now we prescribe
the plumbing family $\{R_t\}$ over the $t$-disc.  Let $D_c=\{|t|<c^4\}$,
$M=R^*_{c^2}\times D_c$ and $\mathcal{V}_c=\{(z,w,t)\mid zw=t,\, |z|, |w|<c
\mbox{ and }|t|<c^4\}$.  $M$ and $\mathcal{V}_c$ are complex manifolds with
holomorphic projections to $D_c$.
Consider the holomorphic maps from $M$ to $\mathcal{V}_c$:
$\hat{F}:(q,t)\rightarrow (F(q),t/F(q),t)$ and
$\hat{G}:(q',t)\rightarrow (t/G(q'),G(q'),t)$ . The maps are consistent with
the projections to $D_c$. The identification space $M\cup \mathcal{V}_c
/\{\hat{F},\hat{G}\mbox{ equivalence}\}$ is a degenerating family $\{R_t \}$
with a projection to the disc $D_c$.  By construction the $0$-fibre has a node
with local model $\mathcal{V}_c$.

We are ready to describe a {\em local manifold cover} of the compactified
moduli space $\Mbar$.  For $R$ having nodes, $R_0=R-\{nodes\}$ is a union of
Riemann surfaces with punctures.  The quasiconformal deformation space of
$R_0$, $Def(R_0)$, is the product of the Teichm\"{u}ller spaces of the
components of $R_0$.  As already noted from \cite[Lemma 1.1]{Wlspeclim} for $m
= \dim\,Def(R_0)$ there is a real analytic family of Beltrami differentials
$\hat{\nu}(s)$, $s$ in a neighborhood of the origin in $\mathbb{C}^m$, such
that $s\rightarrow R_s=R^{\hat{\nu}(s)}$ is a coordinate parameterization of a
neighborhood of $R_0$ in $Def(R)$ and the prescribed mappings
$\zeta^{\hat{\nu}(s)}:R_0\rightarrow R^{\hat{\nu}(s)}$ preserve the cusp
coordinates at each puncture.  Further for $R$ with $n$ nodes we now prescribe
the plumbing data $(U_k,V_k,z_k,w_k,t_k),\,k=1,\dots,n$, for
$R^{\hat{\nu}(s)}$, where $z_k$ on $U_k$ and $w_k$ on $V_k$ are cusp
coordinates relative to the $R^{\hat{\nu}(s)}$-hyperbolic metric (the plumbing
data varies with $s$).  The parameter $t_k$ parameterizes {\em opening} the
$k$ th node.  For all $t_k$ suitably small, perform the $n$ prescribed
plumbings to obtain the family $R_{s,t}=R^{\hat{\nu}(s)}_{t_1,\dots,t_n}$.
The tuple $(s,t)=(s_1,\dots,s_m,t_1,\dots,t_n)$ provides real analytic local
coordinates, the {\em hyperbolic metric plumbing coordinates}, for the local
manifold cover of $\Mbar$ at $R$, \cite{Msext,Wlcut} and \cite[Secs. 2.3,
2.4]{Wlhyp}. The coordinates have a special property: for $s$ fixed the
parameterization is holomorphic in $t$.  The property is a basic feature of
the plumbing construction.  The family $R_{s,t}$ parameterizes the small
deformations of the marked noded surface $R$.

The roles of the Fenchel-Nielsen coordinates and the hyperbolic metric
plumbing coordinates can be interchanged.  In particular for the nodes of $R$
given by the $\sigma$-null stratum $\{\alpha_1,\dots,\alpha_n\}$ the above
local manifold cover has topological coordinates
$((\ell_{\beta},\theta_{\beta})_{\beta \notin
\sigma},(\ell_{\alpha}e^{i\theta_{\alpha}})_{\alpha\in\sigma})$.  The
observation can be established by expressing the Fenchel-Nielsen coordinates
solely in terms of geodesic-lengths, and then applying techniques for
theta-series to analyze the differentials of geodesic-lengths.  Upon
interchanging the roles of the coordinates, we obtain a local description of
the bordification  in terms of the $(s,t)$ tuple, \cite{Abdegn, Bersdeg,
Wlhyp}.  At the point $R$ of the $\sigma$-null stratum in $\Tbar$ the local
parameters are $(s,|t|,\arg t)$ with the $arg$ valued in $\mathbb{R}$.

We review the geometry of the local manifold covers.  For a complex manifold
$M$ the complexification $\boldsymbol{T}^{\mathbb{C}} M$ of the
$\mathbb{R}$-tangent bundle  is decomposed into the subspaces of {\em
holomorphic} and {\em antiholomorphic} tangent vectors.   A Hermitian metric
$g$ is prescribed on the holomorphic subspace.  For a general complex
parameterization  $s=u+iv$ the coordinate $\mathbb{R}$-tangents are expressed
as $\frac{\partial}{\partial u}=\frac{\partial}{\partial
s}+\frac{\partial}{\partial\bar{s}}$ and $\frac{\partial}{\partial
v}=i\frac{\partial}{\partial s}-i\frac{\partial}{\partial\bar{s}}$.   For the
$R_{s,t}$ parameterization the $s$-parameters are not holomorphic while for
$s$-parameters fixed the $t$-parameters are holomorphic;
$\{\frac{\partial}{\partial
s_j}+\frac{\partial}{\partial\bar{s}_j},\,i\frac{\partial}{\partial
s_j}-i\frac{\partial}{\partial\bar{s}_j},\,\frac{\partial}{\partial t_k},\,i
\frac{\partial}{\partial t_k}\}$ is a basis over $\mathbb{R}$ for the tangent
space of the local manifold cover.  For a smooth Riemann surface the dual of
the space of holomorphic tangents is the space of quadratic differentials.
The following is now a modification  of Masur's result \cite[Prop.
7.1]{Msext}.

\begin{proposition}
\label{coords}
The hyperbolic metric plumbing coordinates (s,t) are real analytic and for s
fixed the parameterization is holomorphic in t. Provided the modification
$\hat{\nu}$ is small, for a neighborhood of the origin there are families in
$(s,t)$ of regular 2-differentials $\varphi_j$, $\psi_j$, j=1,\dots,m and
$\eta_k$, k=1,\dots,n  such that: \begin{enumerate}
\item For $R_{s,t}$ with $t_k\ne 0$, all k, $\{\varphi_j,\psi_j,\eta_k\,,i\eta_k\}$ forms the dual basis to  
$\{\frac{\partial\hat{\nu}(s)}{\partial s_j}+\frac{\partial\hat{\nu}(s)}{\partial \bar s_j},i\frac{\partial\hat{\nu}(s)}{\partial s_j}-i\frac{\partial\hat{\nu}(s)}{\partial \bar s_j}, \frac{\partial}{\partial t_k},i\frac{\partial}{\partial t_k}\}$ over $\mathbb{R}$.
\item For $R_{s,t}$ with $t_k=0$, all k, the $\eta_k$, k=1,\dots,n, are trivial and the $\{\varphi_j,\psi_j\}$ span the dual of the holomorphic subspace $\boldsymbol{T}Def(R_0)$.
\end{enumerate}
\end{proposition}

\prf  The situation compares to that considered by Masur.  The new element:
the variation of the plumbing data is prescribed by a Schiffer variation for a
gluing-function real analytically depending on the parameter $s$, \cite[pg.
410]{Wlcut}.  As already noted  for $s$ fixed, plumbing produces a holomorphic
family.  Following Masur the families of regular $2$-differentials
$\{\varphi_j,\psi_j,\eta_k\}$ are obtained by starting with a local frame
$\mathcal{F}$ of regular $2$-differentials and prescribing the  pairings with
$\{\frac{\partial\hat{\nu}}{\partial s_j},\frac{\partial\hat{\nu}}{\partial
\bar s_j},\frac{\partial}{\partial t_k}\}$, \cite[Sec. 5 and Prop.
7.1]{Msext}. At an initial point the  basis is simply given by a linear
transformation of the frame $\mathcal{F}$.  The prescribed basis will then
exist in a neighborhood provided the pairings are continuous.  We first
consider the pairings with $\frac{\partial}{\partial t_k}$. From  (\ref{belt})
we have the Beltrami differential for the  pairing with $\frac{\partial}{\partial t_k}$,
$k=1,\dots,n$.  In particular for a plumbing collar of $R_{s,t}$ let $z$ (or
$w$) be the coordinate of the plumbing.  A quadratic differential $\varphi$ on
$R_{s,t}$ can be factored on the collar into a product of
$\bigl(\frac{dz}{z}\bigr)^2$ and a function holomorphic in $z$.  We write
$\mathcal{C}_k(\varphi)$ for the constant coefficient of the Laurent expansion
of the function factor.  From (\ref{pair}) the pairing with
$\frac{\partial}{\partial t_k}$ is the linear functional
$-\frac{\pi}{t_k}\mathcal{C}_k$.  From Masur's considerations \cite[Sec. 5,
esp. 5.4, 5.5]{Msext} the pairing of $\frac{\partial}{\partial t_k}$ with the local frame $\mathcal F$ is continuous,
and there are  regular $2$-differentials $\{\varphi_j,\psi_j,\eta_k^*\}$ with:
$\mathcal C_{\ell}(\varphi_j)=\mathcal C_{\ell}(\psi_j)=0$, $j=1,\dots,m$;
$\mathcal{C}_{\ell}(\eta^*_k)=\delta_{k\ell}$, $k,\ell=1,\dots,n$.  The
$2$-differentials $\eta_k=\frac{-t_k}{\pi}\eta^*_k,\,k=1,\dots,n$ have the
desired pairings with $\frac{\partial}{\partial t_k}$ .  The final matter is
to note that the pairings of $\{\varphi_j,\psi_j,\eta_k^*\}$  
with $\{\frac{\partial\hat{\nu}}{\partial
s_{\ell}},\frac{\partial\hat{\nu}}{\partial \bar s_{\ell}}\}$ are indeed
continuous in $(s,t)$. By construction the differential $\hat{\nu}(s)$  is
supported in the complement of the plumbing collars, \cite[Lemma
1.1]{Wlspeclim}.  On the support of $\hat{\nu}(s)$ the $2$-differentials are
real analytic in $(s,t)$.  The pairings are continuous and even real analytic.
The proof is complete.

We now note two general matters: the role of the coefficient functional
$\mathcal{C}$, and the approximation of the hyperbolic metric.  As above, for
$z$ a plumbing collar coordinate for $R_{s,t}$, a quadratic differential
$\psi$  can be factored on the collar as the product of
$\bigl(\frac{dz}{z}\bigr)^2$ and a holomorphic function.  $\mathcal{C}(\psi)$
denotes the constant coefficient of the Laurent expansion of the function.
The surface $R_{s,t}$ is constructed by plumbing $\bigl(R_s\bigr)^*_{c^2}$
with the $R_s$-hyperbolic cusp coordinates.  $R_{s,t}$ is the disjoint union
of $\bigl(R_s\bigr)^*_c$, $R_s$ with the cuspidal discs $|z_*|,|w_*|<c$
removed, and the annulus $\{|t|/c<|z|<c\}$.  An approximate hyperbolic metric
$d\omega^2$ is given by choosing the $R_s$-hyperbolic metric on
$\bigl(R_s\bigr)^*_c$ and $dh_t^2$ on the annulus (see (\ref{hyp})).  The
metric $d\omega^2$ is the {\em model grafting } treated in detail in
\cite[Sec. 3.4.MG]{Wlhyp}; as noted in \cite[pgs. 445, 446]{Wlhyp} for
$dh_{s,t}^2$ the $R_{s,t}$-hyperbolic metric we have that
$\bigl|d\omega^2/dh_{s,t}^2\,-\,1\bigr|$ is $O\bigl(\sum_k (\log
|t_k|)^{-2}\bigr)$.  The approximation $d\omega^2$ will now be substituted for
the construction of \cite[Sec. 6]{Msext} to obtain an improved form of the
original expansion.  The improved approximation of the hyperbolic metric is
the new contribution.  Yamada \cite{Yam} presented a third-order expansion
based on the technical work of Wolf \cite{Wfinf} and Wolf-Wolpert
\cite{WWana}.

\begin{theorem}
\label{wp}
For a noded Riemann surface R the hyperbolic metric plumbing coordinates for
$R_{s,t}$ provide real analytic coordinates for a local manifold cover
neighborhood for $\Mbar$.  The parameterization is holomorphic in $t$ for $s$
fixed. On the local manifold cover the WP metric is formally Hermitian
satisfying: \begin{enumerate}
\item For $t_k=0$, $k=1,\dots,n$, the restriction of the metric is a smooth K\"{a}hler metric, isometric to the WP product metric for a product of Teichm\"{u}ller spaces.
\item For the tangents $\{\frac{\partial}{\partial s_j},\frac{\partial}{\partial \bar s_j},\frac{\partial}{\partial t_k}\}$ and the quantity $\rho = \sum\limits^n_{k=1}(\log |t_k|)^{-2}$ then:
$g_{WP}(\frac{\partial}{\partial t_k},\frac{\partial}{\partial t_k})(s,t)=\dfrac{\pi^3}{|t_k|^2(-\log^3 |t_k|)}\,(1\,+\,O(\rho))$;\\
$g_{WP}(\frac{\partial}{\partial t_k},\frac{\partial}{\partial t_{\ell}})\quad  is\quad  O((|t_kt_{\ell}|\log^3|t_k|\log^3|t_{\ell}|)^{-1})\ for\ k \ne \ell$;\\
and for $\mathfrak{u}= \frac{\partial}{\partial s_j},\,\frac{\partial}{\partial \bar s_j}$;\ $\mathfrak{v}= \frac{\partial}{\partial s_{\ell}}\,,\frac{\partial}{\partial \bar s_{\ell}} $:\\ 
$g_{WP}(\frac{\partial}{\partial t_k},\mathfrak{u})\quad is \quad O((|t_k|(-\log^3|t_k|))^{-1})$ and\\
$g_{WP}(\mathfrak{u},\mathfrak{v})(s,t)=g_{WP}(\mathfrak{u},\mathfrak{v})(s,0)\,(1\,+\,O(\rho))$.
\end{enumerate}
\end{theorem}

\prf   We begin with the expansion of the dual metric for the basis provided
in Proposition \ref{coords}.  The behavior of the $\varphi_j$, $\psi_j$,
$\eta_k$ and their contribution to the Petersson pairing
$\int\alpha\overline{\beta}(d\omega^2)^{-1}$ is straightforward.  On
$\bigl(R_s\bigr)^*_{c^2}$ the quadratic differentials and the approximating
metric are real analytic in $(s,t)$.  The contributions to the pairing are
real analytic and each differential $\eta_k,\,k=1,\dots,n$, contributes a
factor of $t_k$.  On the plumbing collars $\{|t|/c<|z|<c\}=\{|t|/c<|w|<c\}$
each quadratic differentials is given as the product of
$\bigl(\frac{dz}{z}\bigr)^2=\bigl(\frac{dw}{w}\bigr)^2$ and a function factor.
We begin with  an elementary calculation

\begin{gather*}
\int_{\{|t|/c<|z|<c\}}\bigl|z^{\alpha}\bigl(\frac{dz}{z}\bigr)^2\bigr|^2\,\bigl(dh^2_t\bigr)^{-1}=\\
\frac{2}{\pi}\int^c_{|t|/c}\bigl(\log |t| \sin \frac{\pi \log r}{\log |t|}\bigr)^2r^{2\alpha}\,d\log r \\
\mbox{where for }\alpha=0,\,\mu=\frac{\log r}{\log |t|}\mbox{ and }\epsilon=\frac{\log c}{\log |t|},\mbox{ then}\\
=\frac{2}{\pi}(-\log^3|t|)\int^{1-\epsilon}_{\epsilon}\sin^2\pi\mu\,d\mu\,=\,\frac{1}{\pi}(-\log^3|t|)\,+\,O(1),\\
\mbox{ and for }\alpha=1, \mbox{ since }|\sin\mu|\le |\mu|, \mbox{ then }\\
=O(1).
\end{gather*}

We are ready to consider the contribution to the Petersson pairing from the
collars.  Consider the contribution for the $\ell^{\mbox{\scriptsize{\,th}}}$
collar.  By construction $\eta_{\ell}$ is the unique quadratic differential
from the dual basis with a nonzero $\mathcal{C}_{\ell}$ evaluation.  In
particular $\mathcal{C}_{\ell}(\eta^*_{\ell})=1$ and the contribution to the
self pairing for $\eta^*_{\ell}$ is $\frac{1}{\pi}(-\log^3|t_{\ell}|)+O(1)$.
In general we note that a quadratic differential on a plumbing collar can be
factored as $\bigl(\frac{dz}{z}\bigr)^2(f_z\,+\,\boldsymbol{c}\,+\,f_w)$ for
$f_z$ holomorphic in $|z|<c,\,f_z(0)=0$; $\boldsymbol{c}$ the
$\mathcal{C}$-evaluation value and $f_w$ holomorphic in $|w|<c,\,f_w(0)=0$.
Furthermore $f_z$, resp. $f_w$, is given as the Cauchy integral of $f$ over
$|z|=c$, resp. $|w|=c$.  Further from the Schwarz Lemma $|f_z|\le
c'|z|\max_{|z|=c}|f|$ with a corresponding bound for $|f_w|$.  The bounds are
combined with the majorant bound $|\sin\mu|\le|\mu|$ to show that: for
$\varphi_j,\,\psi_j,\,\eta^*_k$ on $|z|=c+\epsilon_0$ and $|w|=c+\epsilon_0$
depending analytically on $(s,t)$ their contribution to the Petersson pairing
over the  collar is also analytic in $(s,t)$.

Combining our considerations and noting the approximation of $d\omega^2$ to
the hyperbolic metric  for $R_{s,t}$ we find that

\begin{gather*}
\begin{split}
\bigl<\eta^*_k,\,\eta^*_k\bigr>_{WP}\,&=\,\frac{1}{\pi}(-\log^3|t_k|)\,(1\,+\,O(\sum\limits^n_{\ell=1}(\log|t_{\ell}|)^{-2})),\\
\bigl<\eta^*_k,\,\eta^*_{\ell}\bigr>_{WP}\,&=\,O(1)\mbox{  for  }k\ne\ell,\\
\mbox{and for }\mathfrak{a}=\varphi_j,\psi_j;\, \mathfrak{b}=\varphi_{\ell},\psi_{\ell}:\\
\bigl<\mathfrak{a},\,\eta^*_k\bigr>_{WP}\,&=\,O(1)\mbox{  and}\\
\bigl<\mathfrak{a},\mathfrak{b}\bigr>_{WP}(s,t)\,&=\,\bigl<\mathfrak{a},\mathfrak{b}\bigr>_{WP}(s,0)\,(1\,+\,O(\sum\limits^n_{k=1}(\log|t_k|)^{-2})).
\end{split}
\end{gather*}
The desired expansion now follows from the following Proposition and the
relations $\eta_k=-\frac{t_k}{\pi}\eta_k^*$.  The proof is complete.

For $\mathcal{A}$ a symmetric $m+n\times m+n$ matrix

\begin{equation*}
\begin{aligned}
\begin{pmatrix}
\lambda_1	&	\dots		&a_{\ell j}	&\dots	\\
\vdots	&	\lambda_k	&\vdots	&\vdots	\\
a_{j\ell}	&	\dots		&\ddots	&\dots	\\
\dots		&	\dots		&\dots	&	B
\end{pmatrix}
&
\begin{aligned}
&\quad\mbox{with }\lambda_k,\ 1\le k \le n; \\ &\quad a _{j\ell},\ 1\le j\le m+n,\ j\ne \ell,\ 1\le \ell\le n \mbox{ and} \\ &\quad B\,=\,(b_{j\ell})\mbox{ a symmetric }m\times m \mbox{ matrix},
\end{aligned}
\end{aligned}
\end{equation*}
we consider the situation that $\lambda_1,\dots,\lambda_n$ are large compared to the $a_{j\ell}$ and $b_{j\ell}$.

\begin{proposition}
For\ $\det B\ne0$, and $\rho = \sum\limits^n_{k=1}\lambda_k^{-1}$ then:
$$\det \mathcal{A}=\det B\prod\limits^n_{k=1}\lambda_k(1+O(\rho))$$
and $\mathcal{A}^{-1}=(\alpha_{j\ell})$ where: for  $1\le k\le n,\ \alpha_{kk}=\lambda_k^{-1}(1+O(\rho))$; for $1\le j<\ell\le n,\ \alpha_{j\ell}$ is $O((\lambda_j\lambda_{\ell})^{-1})$; for $1\le j\le n<\ell\le m+n,\ \alpha_{j\ell}$ is $O(\lambda_j^{-1})$, and for $1\le j,\,\ell\le m,\ \alpha_{j+n\,\ell+n}=b^{j\ell}(1+O(\rho))$.  
The constants for the $O$-terms are bounded in terms of $m+n$, $\det B^{-1}$ and $\max \{|a_{j\ell}|,|b_{j\ell}|\}$.
\end{proposition}

\prf We consider the general formula for the determinant as a sum over the
permutation group and by the cofactor expansion.  First observe that there is
a dichotomy for $m+n$-fold products in the calculation of $\det \mathcal{A}$;
a product either also occurs in the expansion of $\det B\prod_k\lambda_k$, or
has at most $n-1$ factors $\lambda_k, 1\le k\le n$.  Products with less than
$n$ factors are bounded in terms of the cited product and $O(\rho)$.  The
determinant expansion is a  consequence. We continue and apply the analog of
the dichotomy when examining the cofactors of $\mathcal{A}$.  For the cofactor
for $\lambda_{\ell}$ we find the expansion $\det
B\prod_{k\ne\ell}\lambda_k(1+O(\rho))$.  Similarly for the cofactor of
$a_{j\ell}$ we find the $\lambda$-contribution to be $\prod_{k\ne
j,\ell}\lambda_k(1+O(\rho))$ for $j\ne\ell\le n$ and to be  $\prod_{k\ne
j}\lambda_{k}(1+O(\rho))$ for $j\le n<\ell$.  Finally for the cofactor of
$b_{j\ell}$ we find the expansion $b^{j\ell}\det B\prod_k\lambda_k(1+O(\rho))$
in terms of the inverse 
$B^{-1}=(b^{j\ell})$.  The proof is complete.

By way of application we present a normal form for the quadratic form
$dg_{WP}^2$; the result is an immediate consequence of the above Theorem.

\begin{corollary} \label{wpnormal} For the prescribed hyperbolic metric
plumbing coordinates: $$
dg^2_{WP}(s,t)\,=\,\bigl(dg^2_{WP}(s,0)\,+\,\pi^3\sum\limits^n_{k=1}(4dr_k^2+r_k^6d\theta^2_k)\bigr)\,(1+O(\|r\|^3))
$$ for $r_k=(-\log|t_k|)^{-1/2},\ \theta_k=\arg t_k$ and $r=(r_1,\dots,r_n)$.\
\end{corollary}

The result provides a local expansion of the WP metric about the
compactification divisor $\mathcal{D}=\{t_k=0\}$. To the third order of
approximation the WP metric is formally a product.  As we will note below, a
second-order approximation is already special.  As already noted, the
bordification $\Tbar$  has a local description in terms of the parameters
$(s,|t|,\arg t)$ or equivalently in terms of $(s,(-\log |t|)^{-1/2},\arg t)$.
The above result provides the associated WP expansion.

An almost-product Riemannian metric with remainder bounded by the displacement
from a submanifold is very special.  We note the situation as motivation for
the results of Section \ref{lengthmin}; the following considerations do not
apply since $4dr^2+r^6d\theta^2$ is not a Riemannian metric.   Consider a
product $\mathbb{R}^m\times\mathbb{R}^n$ with Euclidean coordinates $x$ for
$\mathbb{R}^m$ and $y$ for $\mathbb{R}^n$.  Consider that in a neighborhood of
the origin a metric has the expansion

$$ dg^2\,=\,dg_x^2\,+\,dg^2_y\,+\,O_{C^1}(\|y\|^2) $$
with $dg_x^2$, resp. $dg_y^2$, a $C^1$-metric for $\mathbb{R}^m$, resp. for
$\mathbb{R}^n$, and the remainder a $C^1$-symmetric tensor as indicated.  The
expansion provides that the {\em second fundamental form} of the $x$-axes,
$\mathbb{R}^m\times \{0\}$, vanishes identically \cite[pgs. 62, 100]{BO'N}.
In this case the $x$-axes is a {\em totally geodesic} submanifold: a geodesic
initially tangent to the $x$-axes is contained in the $x$-axes \cite[pg.
104]{BO'N}.  The  expansion also provides that for the $x$-axes the {\em
normal connection} and the {\em normal curvature} vanish identically
\cite[pgs. 114, 115]{BO'N}.

\section{Length-minimizing curves on Teichm\"{u}ller space}

We begin by developing basic facts about the behavior of WP geodesics on
Teichm\"{u}ller space.  Although Teichm\"{u}ller space is topologically a
cell, the behavior of geodesics is not a consequence of general results
\cite{KN}, since the WP metric is not complete.  For instance the Hopf-Rinow
theorem cannot be directly applied to obtain length-minimizing curves
\cite{CE,KN,BO'N}, and it is necessary to show that distance is measured along
geodesics.  We proceed though by applying our paradigm: {\em a filling
geodesic-length sum behaves qualitatively as the distance from a point for a
complete metric.}  In the following we combine the paradigm and modifications
of the standard arguments to find the basic behavior of geodesics.

\begin{theorem} \label{WPexp} The WP exponential map from a base point is a
diffeomorphism from its open domain onto the Teichm\"{u}ller space.
\end{theorem}

\begin{corollary} \label{unqgeo} Teichm\"{u}ller space is a unique geodesic
space.  Each WP geodesic segment is the unique length-minimizing rectifiable
curve connecting its endpoints.  \end{corollary}

\noindent{\bf Proof of Corollary.}\  Let $\gamma$ be the WP geodesic
connecting a pair of points $p$ and $q$ in the Teichm\"{u}ller space.  For a
filling geodesic-length function $\mathcal{L}$, choose $c>0$, such that
$\gamma\subset S_c=\{\mathcal L<c\}$, \cite{Wlnielsen}.  Consider $G$ the set
of all rectifiable curves connecting $p$ and $q$, contained in $\overline
{S_c}$, and each with length at most $d(p,q)+ 1$.  Provided $G$ is nonempty
and the elements of $G$ are parameterized proportional to arc-length on the
interval $[0,1]$, then $G$ constitutes an equicontinuous family of maps.  In
particular for $\beta\in G$ and $t,t'\in [0,1]$ by the proportional
parameterization it follows that

$$ |t-t'|=\frac{L(\beta [t,t'])}{L(\beta)}\ge
\frac{d(\beta(t),\beta(t'))}{d(p,q)+1}.  $$
From the Arzel\`{a}-Ascoli Lemma \cite[pg. 36]{BH} there exists a rectifiable
length-minimizing (amongst elements of $G$) curve $\beta_0$ connecting $p$ and
$q$ contained in $\overline{S_c}$.

We  consider the behavior of a  rectifiable length-minimizing (amongst
elements of $G$) curve $\beta_0$ passing through an arbitrary point $r\in
\overline{S_c}$ ($r$ could lie on $\partial\overline{S_c}$).  Since
$\overline{S_c}$ is compact, there is a positive  $\epsilon$ such that  WP
geodesics are uniquely length-minimizing in an $\epsilon$-neighborhood  of
each point of $\overline{S_c}$.  Since $\overline{S_c}$ is WP convex it
follows for $r', r''$ on the chosen curve, close to $r$, with $r'$ {\em
before} $r$ and $r''$ {\em after} $r$, that the segments $\widehat{r'r}$ and
$\widehat{rr''}$ are necessarily WP geodesics.  It further follows that
$\widehat{r'r''}$ is a WP geodesic, since  the segment is locally
length-minimizing at $r$.  By convexity of the geodesic-length function
$\mathcal L$, its value at $r'$ or $r''$ is greater than its value at $r$.
Since $r',r''\in \overline{S_c}$ it follows that $r\in S_c$.  It now follows
that a rectifiable length-minimizing (amongst elements of $G$)  curve
$\beta_0$ is a WP geodesic entirely contained in $S_c$.

In general given $\epsilon$, $0<\epsilon<1$, there exists a curve $\beta'$
connecting $p$ and $q$ in the Teichm\"{u}ller space such that
$L(\beta')<d(p,q)+\epsilon$.  For $c'$ large, $\beta' \subset \overline{
S_{c'}}$ and thus the corresponding family of maps $G$ is non empty. The
length of $\beta'$ bounds the length of a $S_{c'}$-length-minimizing curve
$\beta_0$ connecting $p$ and $q$: in particular  $L(\beta')\ge L(\beta_0)$.
From the above paragraph and the Theorem, the unique geodesic connecting $p$
and $q$ is $\beta_0=\gamma$.   The inequalities now provide that
$L(\gamma)<d(p,q)+\epsilon$.  The proof is complete.  \\

\noindent{\bf Proof of Theorem.}\  First we note that the domain of the
exponential map is an open set.  Given a geodesic $\gamma$ connecting a pair
of points, select $c>0$ such that $\gamma \subset S_c$.  Since $S_c$ is open
the points in  neighborhoods of the $\gamma$-endpoints are also connected by
WP geodesics.  In particular the domain of the exponential map is open.

We next note that the exponential map is a local diffeomorphism \cite{KN}.
Further note that a germ of the inverse is determined by its value at a single
point.  We now consider the continuation of a given germ $\boldsymbol \iota$,
with the exponential map based at $p$ and the germ given at $q\in \mathcal T$.
We consider the continuation of $\boldsymbol \iota$ along $\alpha$, a curve
with initial point $q$.  We argue that the continuation set is closed.

Choose a filling geodesic-length function $\mathcal L$ and value $c$ such that
$p, \alpha \subset \overline{S_c}$.  First we observe that each WP geodesic
connecting $p$ and a point of $\alpha$ is contained in $\overline{S_c}$.  This
follows since the values of $\mathcal L$ at the endpoints are bounded by $c$
and $\mathcal L$ is WP convex. Since $\overline{S_c}$ is compact there is an
overall length bound for the WP geodesics contained in $\overline{S_c}$.  As
noted in Bridson-Haefliger a length-bounded family of geodesics is given by an
equicontinuous family of maps, \cite[pg. 36]{BH}.  By the Arzel\`{a}-Ascoli
Lemma it follows that a sequence of WP geodesics contained in $\overline{S_c}$
has a subsequence converging to a geodesic contained in $\overline{S_c}$.

Consider now that the germ $\boldsymbol \iota$ can be continued to a sequence
of points $\{q_n\}$ along $\alpha$.  In particular WP geodesics
$\widehat{pq_n}$ are determined.  A subsequence (same notation)
$\widehat{pq_n}$ converges to $\widehat{pq'}$.  The WP geodesic
$\widehat{pq'}$ determines a germ of the inverse of the exponential map;  the
germ gives {\em exponential inverses} for the WP geodesics $\widehat{pq_n}$.
The germ is the continuation of $\boldsymbol \iota$; the continuation set is
closed.  The continuation set is necessarily open; $\boldsymbol \iota$ can be
continued along every curve.  On considering homotopies it is established that
the continuation to the endpoint of $\alpha$ is path independent.  Finally
since the Teichm\"{u}ller space is simply connected the continuations
determine a global inverse for the exponential map.  The proof is complete.

We are also interested in understanding the WP {\em join} of two sets, and in
particular the distance between points on a pair of geodesics. For the WP
inner product consider the Levi-Civita connection $\nabla$  satisfying for
vector fields $X$, $Y$ and $W$ the relations,\cite{CE,KN, BO'N},

\begin{gather*}
X\bigl<Y,W\bigr>=\bigl<\nabla_XY,W\bigr>+\bigl<Y,\nabla_XW\bigr>\\
\nabla_XY-\nabla_YX=[X,Y].  \end{gather*}
Further consider the  curvature tensor

$$ R(X,Y)Z=\nabla_X\nabla_YZ-\nabla_Y\nabla_XZ-\nabla_{[X,Y]}Z.  $$
A {\em variation of geodesics} is a smooth map $\beta(t,s)$ from
$[t_0,t_1]\times (-\epsilon, \epsilon)$ to $\mathcal T$ such that for each
$s'$, $\beta(t,s')$ is a WP geodesic.  For the vector fields
$T=d\beta(\frac{\partial}{\partial t})$ and $V=d\beta(\frac{\partial}{\partial
s})$, the first variation of the geodesic $V$ satisfies the Jacobi equation

$$ \nabla_T\nabla_TV=R(T,V)T; $$
solutions are Jacobi fields, \cite{CE}.  The Jacobi equation is a linear
second-order system of ordinary differential equations for vector fields along
$\beta(t,s')$.  The space of solutions has dimension $2\dim \mathcal T$; a
solution is uniquely prescribed by its initial value and its initial
derivative.  Furthermore since $\mathcal T$ has negative curvature there are
no conjugate points along a geodesic and the linear map $(V\vert
_{t=t_0},\nabla_TV\vert_{t=t_0})$ to $(V\vert_{t=t_0},V\vert_{t=t_1})$ is an
isomorphism \cite[pg. 19]{CE}.  Solutions are  uniquely prescribed by their
values at the endpoints.  This property is needed for the understanding of the
join of sets.  In particular Jacobi fields provide a mechanism for analyzing
the exponential map.

The WP exponential map $(p,v)\stackrel{\boldsymbol{e}}{\rightarrow} exp_pv$
has domain $\boldsymbol{\mathcal D}$ the open set $\{(p,exp_p^{-1}(\mathcal
T))\}\subset \mathbf T\mathcal T$ of the tangent bundle.  We are ready to
consider the behavior of geodesics.

\begin{proposition} \label{discon} The WP exponential map
$(p,v)\stackrel{\boldsymbol{e}}{\rightarrow} (p,exp_pv)$ is a diffeomorphism
from $\boldsymbol{\mathcal D}\subset \mathbf T\mathcal T$ to $\mathcal T\times
\mathcal T$.  For a pair of disjoint WP geodesics parameterized proportional
to arc-length, the WP distance between corresponding points is a strictly
convex function.  \end{proposition}

\prf The map $\boldsymbol{e}$ is smooth with differential
$d\boldsymbol{e}=(id,d\, exp)$.  As already noted since the WP metric has
negative curvature $d\, exp$ has maximal rank and thus $\boldsymbol{e}$ is a
local diffeomorphism of $\mathcal D$ to $\mathcal T\times \mathcal T$.  A
consequence of Corollary \ref{unqgeo} is that $\boldsymbol{e}$ is a global
diffeomorphism.

We are ready to consider the distance between corresponding points of a pair
of disjoint WP geodesics.  From the above result a one-parameter variation of
geodesics is determined $\beta(t,s)$, $(t,s)\in [t_0,t_1]\times [s_0,s_1]$.
For a value $s'\in [s_0,s_1]$ we write $T$ for the tangent field of
$\beta(t,s')$ and $V$ for its variation field; we assume $\|T\|=1$.  The
second variation in $s$ at $s'$ of the length of $\beta(t,s)$ is given by the
classical formula \cite[(1.14)]{CE}

\begin{multline} \label{secondvar}
\bigl<\nabla_VV,T\bigr>\bigr{\vert}^{t_1}_{t_0}+\\
\int^{t_1}_{t_0}\bigl<\nabla_TV,\nabla_TV\bigr>-\bigl<\nabla_TV,T\bigr>\bigl<\nabla_TV,T\bigr>-\bigl<R(V,T)T,V\bigr>
dt.  \end{multline}
Observations are in order.  First by hypothesis the curves $\beta(t_0,s)$ and
$\beta(t_1,s)$ are geodesics with constant speed parameterization; the
acceleration $\nabla_VV$ vanishes at $t_0$ and $t_1$.  Second the first two
terms of the integrand combine to give the length-squared of the projection
$\nabla_TV$ onto the normal space of $T$.  And the third term of the integrand
is strictly positive given strictly negative  curvature \cite{CE}.  In summary
the distance is a strictly convex function.  The proof is complete.

We are ready to show that $\overline{\mathcal T}$ is a geodesic space.  For
points $p$ and $q$ of the completion let $\{p_n\}$ and $\{q_n\}$ be sequences
from $\mathcal T$ converging to $p$, resp. to $q$.  Note for the distance we
have $d(p,q)=\lim_nd(p_n,q_n)$.  Consider the sequence of curves
$\gamma_n=\widehat{p_nq_n}$ of $\mathcal T$ parameterized proportional to
arc-length by the unit-interval.  Since the sequences $\{p_n\}$ and $\{q_n\}$
are Cauchy it follows from Proposition \ref{discon} that for each $t\in[0,1]$
the sequence $\{\gamma_n(t)\}$ is also Cauchy ({\em without} passing to a
subsequence). The sequence $\{\gamma_n\}$ prescribes a function $\gamma$ with
domain the unit-interval and values in $\overline{\mathcal T}$.  Furthermore
since the $\gamma_n$ are distance proportional-parameterized, for
$t,t'\in[0,1]$ then 

$$ \frac{d(\gamma_n(t),\gamma_n(t'))}{d(p_n,q_n)}=|t-t'|.  $$
It follows that $d(\gamma(t),\gamma(t'))=|t-t'|d(p,q)$; $\gamma$ is a
continuous function, in particular a geodesic.  We summarize the
considerations with the following.

\begin{proposition} \label{Tcomp} The completion $\overline{\mathcal T}$ is a
geodesic space.  \end{proposition}

\section{Length-minimizing curves on the completion} \label{lengthmin}

For the  complex of curves $C(F)$ a $k$-simplex is a set of $k+1$ free
homotopy classes of nontrivial, nonperipheral, mutually disjoint simple closed
curves for the reference surface $F$.  A simplex $\sigma$ {\em precedes} a
simplex $\sigma'$ provided $\sigma \subseteq\sigma'$; preceding is a partial
ordering.  With the convention that the {\em $-1$-simplex} is the null set,
there is a natural function $\Lambda$ from the completion $\Tbar$ to the
complex $C(F)\cup\{\emptyset\}$ determined by the classes of the nodes. For a
marked noded Riemann surface $(R,f)$ with $f:F\rightarrow R$, the labeling
$\Lambda((R,f))$ is the simplex of free homotopy classes on $F$ mapped to the
nodes on $R$.  The level sets of $\Lambda$ are the strata of $\Tbar$. We write
$\mathcal S(\sigma)$ for the stratum determined by the simplex $\sigma$.
The stratum for a $k$-simplex has complex dimension $3g(F)-3-k$.  

We now consider first properties of length-minimizing curves on $\Tbar$.  We
are able to make the analysis without first establishing that a
length-minimizing curve is a limit of WP geodesics.  In this section we build
on the following result and present an alternative approach to the basic
observation of S. Yamada \cite{Yam} that except possibly for its endpoints, a
length-minimizing curve is contained in a single stratum of $\Tbar$.

\begin{proposition} \label{finite} For a length-minimizing curve $\gamma$ on
$\Tbar$ the composition $\Lambda\circ\gamma$ has a left and right limit at
each point.  The composition is continuous at a point where the left and right
limits agree.  \end{proposition}

\prf First observe that only a finite number of simplices precede a given
simplex.  There is a continuous analog for strata: in a suitable neighborhood
of a point of $\Tbar$ there are only a finite number of strata, and each
precedes or coincides with the stratum of the point.  If a left or right limit
fails to exist for $\Lambda\circ\gamma$ at $t_0$, then there is a monotonic
convergent sequence of parameter values $\{t_n\},\ t_n\rightarrow t_0$ with
$\Lambda\circ\gamma$ having value $\sigma$ on $\{t_{2n}\}$ and a different
value $\tau$ on $\{t_{2n+1}\}$.  We may choose that $\sigma$ precedes $\tau$
and further that $\sigma$, resp. $\tau$,  is a maximal, resp. minimal, such
value.   Maximal connected segments of $\gamma$ contained in the stratum of
$\sigma$ are determined by the positivity of the geodesic-length functions of
the classes in $\tau-\sigma$.  In particular each maximal segment is
parameterized by an open parameter interval;  by Corollary \ref{unqgeo} each
(closed subsegment of each) maximal segment is length-minimizing.  

Consider two points $\gamma(t_{2n})$ and $\gamma(t_{2n+2})$  on different
maximal segments.  By Corollary \ref{unqgeo} there exists a WP geodesic
$\beta$ contained in the stratum of $\sigma$ connecting $\gamma(t_{2n})$ and
$\gamma(t_{2n+2})$.  We now compare the segments
$\gamma\vert_{[t_{2n},\,t_{2n+2}]}$ and $\beta$.  Assume each segment is
parameterized by arc-length; the segments necessarily have the same length.
On the stratum $\sigma$ the curves $\beta$ and the maximal segment of $\gamma$
at $t_{2n}$ are solutions of an ordinary differential equation.  If the
initial tangents of $\beta$ and $\gamma\vert_{[t_{2n},\,t_{2n+2}]}$ coincide,
then by the uniqueness of solutions and the maximality, the segments must
coincide for the length of $\beta$.  The coinciding contradicts
$\Lambda\circ\gamma$ having different values at $t_{2n}$ and $t_{2n+1}$.  The
alternative is that the initial (unit) tangents of $\beta$ and
$\gamma\vert_{[t_{2n},\,t_{2n+2}]}$ differ.  In this case $\gamma$ can be
modified by first substituting the segment $\beta$ for the parameter interval
$[t_{2n},t_{2n+2}]$ and then {\em smoothing the corner} (inside the stratum)
at $\gamma(t_{2n})$, to obtain a new curve $\widetilde{\gamma}$ of strictly
smaller length, again a contradiction. A sequence $\{t_n\}$ as described
cannot exist.  In summary the composition $\Lambda\circ\gamma$ is locally
constant to the left and right of each point of its domain.

Finally if the left and right limits have a common value at $t_0$ then either
$\Lambda\circ\gamma(t_0)$ also has the common value, or the common value
precedes $\Lambda\circ\gamma(t_0)$.  In the second instance we can again
construct a modification $\widetilde{\gamma}$ of strictly smaller length.  The
proof is complete.

We are interested in a class of singular {\em metrics} that model the WP
metric in a neighborhood of a point on the compactification divisor
$\mathcal{D}\subset\overline{\mathcal{M}}$.  Consider now the product
$(\mathbb{R}^2)^{m+n}$ with Euclidean coordinates $(x,y)$ for $x$ the
$2m$-tuple with Euclidean metric $dx^2$ and $y$ the $2n$-tuple with Euclidean
metric $dy^2$.  We refer to $\mathbb{R}^{2m}\times\{0\}$, resp. to
$\{0\}\times\mathbb{R}^{2n}$, as the $x$-axes, resp. the $y$-axes.  Here the
$x$-axes represent coordinates on a  stratum of dimension $2m$ and codimension
$2n$, while the $y$-axes represent the parameters which open nodes.  We write
$(r_j,\theta_j)$ for the polar coordinates for the $2$-plane
$(y_{2j-1},y_{2j})$ and $(r,\theta)$ for the product polar coordinates for the
$2n$-tuple of  $y$-coordinates.  We consider the singular metric

$$ \sum\limits^n_{j=1}\,4\,dr_j^2+r_j^6\,d\theta_j^2 $$
for the $y$-axes which we simply abbreviate as $dr^2+r^6d\theta^2$.

\begin{definition} \label{cuspmetric} A continuous symmetric $2$-tensor $ds^2$
is a product cuspidal metric for a neighborhood of the origin in
$(\mathbb{R}^2)^{m+n}$ provided: \begin{enumerate} \item  $ds^2$ is a smooth
Riemannian metric on $\bigcap_{j=1}^n\{r_j>0\}$; \item  the restriction of
$ds^2$ to the $x$-axes is a smooth Riemannian metric $ds_x^2$; \item
$ds^2=(d\mu^2+dr^2+r^6d\theta^2)(1+O(\|r\|^2))$\ \ for $d\mu^2$ the pullback
of $ds_x^2$ to $(\mathbb{R}^2)^{m+n}$ by the projection onto the $x$-axes ,and
$\|r\|$ denoting the Euclidean norm of the radius vector for the
$y$-axes.\label{pc3} \end{enumerate} \end{definition}

We are ready to continue our consideration of a length-minimizing curve
$\gamma$ and a point of discontinuity $t_*$ interior to the domain of the
label composition $\Lambda\circ\gamma$.  The first circumstance to consider is
that $\Lambda\circ\gamma$ is continuous from one side, say the right.  In
particular for $\Lambda\circ\gamma$ discontinuous from the left the simplex
$\sigma=\Lambda\circ\gamma(t_*^-)$ strictly precedes the simplex
$\sigma'=\Lambda\circ\gamma(t_*)$.  Since $\gamma$ is length-minimizing, the
curve is a WP geodesic in the stratum $\mathcal S(\sigma')$ for an initial
interval to the right of $t_*$.  Observe by analogy that in the circumstance
that $ds^2,\,d\mu^2$ and $dr^2+r^2d\theta^2$ were smooth, then the $x$-axes
would be totally geodesic (the {\em second fundamental form} would be trivial;
see the discussion after Corollary \ref{wpnormal}) and the suggested
refracting behavior of $\gamma$ would not be possible.  We will now show by a
scaling argument that the behavior is also not possible for a product cuspidal
metric.

We observe that the individual strata of $\Tbar$ branch cover strata in
$\Mbar$ and that certain curves in $\Mbar$ have unique lifts determined by
their initial point.  $\Mbar$ is a $V$-manifold; consider first the
local-manifold cover $\widehat{\mathcal{U}}$ for a neighborhood of the point
given by $\gamma(t_*)$ (a neighborhood of $\gamma(t_*)$ in $\Mbar$ is
described as $\widehat{\mathcal{U}}/Aut(\gamma(t_*))$).  From Section
\ref{prelim} the preimage of $\widehat{\mathcal{U}}$ in $\Tbar$ is a disjoint
union of sets, including a neighborhood $\mathcal{U}$ of $\gamma(t_*)$.  The
local stratum $\sigma\cap\mathcal{U}\subset\Tbar$ is a covering of its image
$\widehat{\sigma}\subset\widehat{\mathcal{U}}$, with covering group the
lattice of Dehn twists for the set of loops $\sigma'-\sigma$.  The local
stratum $\sigma'\cap\mathcal{U}\subset\Tbar$ coincides with its local
projection $\widehat{\sigma'}$ to $\widehat{\mathcal{U}}$.  In the following
paragraphs we will use the simple observation that: a rectifiable curve in
$\widehat{\mathcal{U}}$ with first segment in $\widehat{\sigma}$ and second
segment in $\widehat{\sigma'}$ has WP isometric lifts to $\mathcal{U}$, each
uniquely determined by prescribing an initial point.  We can choose
coordinates so $\widehat{\mathcal{U}}\cap\widehat{\sigma}$ is given by a
neighborhood of the origin in $(\mathbb{R}^2)^{m+n}$ and
$\widehat{\mathcal{U}}\cap\widehat{\sigma'}$ is given by a neighborhood of the
origin in the $x$-axes $\mathbb{R}^{2m}\times\{0\}$.

We study the WP length of an oriented curve $\gamma$ having first segment off
the $x$-axes and second segment a geodesic in the $x$-axes.  Let $\mathbf{o}$
be the first contact point of the curve with the $x$-axes.  Consider the
Euclidean ball of radius $\delta$ about $\mathbf{o}$. From Corollary
\ref{wpnormal} we can choose $\delta$ small for the metric to have the
coordinate description of a product cuspidal metric in the ball. Along
$\gamma$ let $a$ be the first intersection point of $\gamma$ with the
$\delta$-sphere at $\mathbf{o}$ to the left of $\mathbf{o}$. Along  $\gamma$
let $b$ be the first intersection point of $\gamma$ with the $\delta$-sphere
at $\mathbf{o}$ to the  right of $\mathbf{o}$. We investigate the lengths of
curves from $a$ to $\mathbf{o}$ to $b$ as $\delta$ varies.  We use the
coordinates $(x,y)$ for the following constructions.  Let $a_x$, resp. $a_y$,
be the Euclidean projection of $a$ to the $x$, resp. the $y$, axes.  Let
$\beta$ be the unit-speed $d\mu^2$ geodesic in the $x$-axes from $a_x$ to $b$.
For the same arc-length parameter let $\widehat{\beta}$ be the curve from $a$
to $b$ whose Euclidean projection to the $y$-axes is a constant speed radial
line.  On $\widehat{\beta}$ the tensor $dr^2+r^6d\theta^2$ restricts  to
$dr^2$ and in particular 

$$
\int_{\widehat{\beta}}(ds^2)^{1/2}=\int_{\widehat{\beta}}(1+O(\|y\|^2))(d\mu^2+dr^2)^{1/2}.
$$
Since the length of $\widehat{\beta}$ is bounded in terms of $\delta$ and the
Euclidean height of $\widehat{\beta}$ is bounded by $\|a_y\|$ it follows that
the length of $\widehat{\beta}$ is given as

$$ \int_{\widehat{\beta}}(d\mu^2+dr^2)^{1/2}+O(\delta\|a_y\|^2) $$
for the Euclidean norm of $a_y$.  The integral immediately evaluates to
$(\trip \beta\trip^2+\|a_y\|^2)^{1/2}$ for $\trip\ \trip$ denoting the
$d\mu^2$ length.  In summary the length of $\widehat{\beta}$ is bounded above
by $(\trip \beta \trip^2+\|a_y\|^2)^{1/2}+O(\delta\|a_y\|^2)$.  

We next consider a lower bound for the length of the segment of $\gamma$ from
$a$ to $b$.  The first expansion is provided similar to the above
consideration.  Select a subsequence of values $\delta'$ tending to zero such
that for each $\delta'$ the maximum of the $y$-height $\|y\|$ on the $\gamma$
segment from $a$ to $\mathbf{o}$ actually occurs at the initial point $a$. For
the subsequence the length of $\gamma$ is

$$ \int_{\gamma}(d\mu^2+ dr^2+r^6d\theta^2)^{1/2}+O(\delta'\|a_y\|^2).  $$
The metric $d\xi^2=(d\mu^2+ dr^2+r^6d\theta^2)$ is a product and the local
behavior of its geodesics is understood. Let $\widehat{\gamma}$ be a
comparison curve (a $d\xi^2$ piecewise geodesic) with segments from $a$ to
$\mathbf{o}$ and from $\mathbf{o}$ to $b$.  The projection of the first
segment of $\widehat{\gamma}$ to the $x$-axes is the $d\mu^2$ geodesic from
$a_x$ to $\mathbf{o}$. The projection of the first segment of
$\widehat{\gamma}$ to the  $y$-axes is  the constant speed radial line from
$a_y$ to $\mathbf{o}$.  The second segment of $\widehat{\gamma}$ is the
geodesic from $\mathbf{o}$ to $b$.  The integral of $\gamma$ is minorized by
the integral of $\widehat{\gamma}$ and consequently the length of $\gamma$ is
minorized by

$$ (\trip a_x\trip^2+\|a_y\|^2)^{1/2}+\trip b\trip^2+O(\delta'\|a_y\|^2) $$ 
where we have written $\trip x_0\trip$ for the $d\mu^2$-distance from
$\mathbf{o}$ to the point $x_0$ of the $x$-axes.

We are prepared to analyze the length of $\gamma$ in a small neighborhood of the point
$\mathbf{o}$.

\begin{proposition} 
\label{refract} A curve having first segment off the
$x$-axes and second segment a geodesic in the $x$-axes is not $ds^2$
length-minimizing between its endpoints.  There is a shorter curve of the same
description.  
\end{proposition}

\prf We first consider the rescaling limit of a neighborhood of $\mathbf{o}$
with the substitution $\delta u=x, \delta v=r$ and $d\eta^2=\delta^{-2}ds^2$.
The curves $\widehat{\beta}$ and $\widehat{\gamma}$ considered above have
radial lines as their projections to the $y$-axes; it suffices for length
considerations to consider the projection of the $y$-axes to its radial
component $r$.  The rescaling limit of $(d\mu^2+dr^2)(1+O(\|y\|^2))$ is the
Euclidean metric and for a subsequence the points $a$, $b$ limit to points of
the unit sphere (same notation).  The curve $\widehat{\beta}$ limits to the
chordal line connecting $a$ to $b$; the curve $\widehat{\gamma}$ limits to the
segmented curve of line-segments connecting $a$ to $\mathbf{o}$ and
$\mathbf{o}$ to $b$.  If $a$ is not antipodal to $b$ then (on the subsequence)
$\widehat{\beta}$ is strictly shorter than $\widehat{\gamma}$. On a
neighborhood of $\mathbf{o}$, $\gamma$ is now modified by substituting a
segment of $\widehat{\beta}$ to obtain a strictly shorter curve, a desired
conclusion.

It remains to consider that the rescaling limit as $\delta$ tends to zero of
$a$ is the antipode to $b$.  In this circumstance we have that $\trip
a_x\trip$ is comparable to $\delta$, $\trip \beta\trip$ is comparable to
$2\delta$ and $\|a_y\|$ by hypothesis is $o(\delta)$.  Pick $\epsilon <1$ such
that $\epsilon^2\trip \beta\trip > \trip a_x \trip$ for all small $\delta$.
Now from the preliminary considerations for small $\delta$ the length of
$\widehat{\beta}$ is bounded above by

$$ (\trip \beta\trip^2+\|a_y\|^2)^{1/2}+O(\delta\|a_y\|^2)\le\trip\beta\trip
+\frac{1}{2\epsilon\trip\beta\trip}\|a_y\|^2+O(\delta\|a_y\|^2) $$
and for a suitable subsequence the length of $\gamma$ is bounded below by

\begin{multline*} (\trip a_x\trip^2+\|a_y\|^2)^{1/2}+\trip
b\trip+O(\delta\|a_y\|^2)\ge \\ \trip a_x\trip+\trip b\trip
+\frac{\epsilon}{2\trip a_x\trip}\|a_y\|^2+O(\delta'\|a_y\|^2) \end{multline*}
for $\|a_y\|(\trip a_x\trip)^{-1}$ sufficiently small, which is ensured for
$\delta$ sufficiently small.  As specified above, $\delta'$ are the special
values for which the maximum of the $y$-height $\|y\|$ on the $\gamma$ sement
from $a$ to $\mathbf{o}$ occurs at the initial point $a$.  Since $\beta$ is a
geodesic we have that $\trip\beta\trip\le\trip a_x\trip +\trip b\trip$.
Observe that the coefficient of the $\| a_y\|^2$-term for $\widehat{\beta}$ is
strictly less than that for $\gamma$.   Since $\| a_y\|^2$ is positive for
$\delta$ positive, it now follows that $\widehat{\beta}$ is strictly shorter
than $\gamma$ and in particular that $\gamma$ is not length-minimizing in a
neighborhood of $\mathbf{o}$.  The proof is complete.

The second circumstance to consider is that for a length-minimizing curve
$\gamma$ there is an interior domain discontinuity point $t_*$ with the label
composition $\Lambda\circ\gamma$ not continuous from either side.  The curve
$\gamma$ connects points in different strata by passing through a higher
codimension stratum.  In particular from Proposition \ref{refract} it follows
that the values of $\Lambda\circ\gamma$ for $t_*^-,t_*$ and $t^+_*$ are all
distinct; furthermore the values for $t_*^-$ and $t_*^+$ strictly precede the
value at $t_*$.  We have the coordinate description of the local strata
$\widehat{\mathcal U}\,\cap\,(\mathcal{S}\circ\Lambda\circ\gamma(t_*^-)\,
\cup\,\mathcal{S}\circ\Lambda\circ\gamma(t_*)\,\cup\,\mathcal{S}\circ\Lambda\circ\gamma(t_*^+))$
given in a neighborhood of the origin in $(\mathbb{R}^2)^{m+n}$.  For suitable
$n_-\,+\,n_+=n$ the neighborhood is given as a neighborhood of the origin in
$(\mathbb R^2)^{m+\,n_-+\,n_+}$ with coordinates $(x,y_-,y_+)$.  In a
neighborhood of the origin the three strata are given by germs of the
coordinate axes: $\mathcal S\circ\Lambda\circ\gamma(t_*)$ by the $x$-axes;
$\mathcal S\circ\Lambda\circ\gamma(t_*^-)$ by the $y_-$-axes; and $\mathcal
S\circ\Lambda\circ\gamma(t_*^+)$ by the $y_+$-axes.  Again a rectifiable curve
with the prescribed behavior for $\Lambda\circ\gamma$ has WP isometric lifts
to $\mathcal{U}$, each uniquely determined now by prescribing an initial and
terminal point.

\begin{proposition} \label{corners} A curve having endpoints distinct from the
origin and in distinct coordinate proper subspaces of the $y$-axes and further
having the origin as an intermediate point is not $ds^2$ length-minimizing
between its endpoints.  There is a shorter curve avoiding the origin.
\end{proposition}

\prf  The considerations are simplified since in effect the subspaces
corresponding to $\mathbb{R}^{2n_-}$ and $\mathbb{R}^{2n_+}$ are orthogonal.
Choose $\epsilon > 0$ such that $(1+\epsilon)<(1-\epsilon)\sqrt{2}$; from
Definition \ref{cuspmetric} the restriction of $ds^2$ to the $y$-axes is
estimated above, resp. below, by the $(1+\epsilon)$, resp. $(1-\epsilon)$,
multiple of $(dr^2+r^6d\theta^2)$.  For $a^-$ in  $\mathbb{R}^{2n_-}-\{0\}$
and $a^+$ in  $\mathbb{R}^{2n_+}-\{0\}$ the $(dr^2+r^6d\theta^2)$ geodesic in
$\mathbb R^{2n}$ connecting $a^-$ to $a^+$ and the piecewise geodesic
connecting $a^-$ to the origin and then to $a^+$ are Euclidean line segments.
The line connecting $a^-$ to $a^+$ has $ds^2$ length at most
$(1+\epsilon)(\|a^-\|^2+\|a^+\|^2)^{1/2}$.  The line segments connecting $a^-$
to the origin to $a^+$ have length at least $(1-\epsilon)(\|a^-\|+\|a^+\|)$.

For an oriented curve $\gamma$ with the prescribed strata behavior consider a
Euclidean radius $\delta$ sphere at the origin and let $a^-$, resp. $a^+$, be
the first intersection point along $\gamma$ to the left, resp. right, of the
origin.  Since the radial component of the metric is comparable to the
Euclidean metric, the maximum value of $\|r\|$ along the segments of $\gamma$
is comparable to $\delta$.  Apply the above estimate for $\delta$ small to
obtain the desired conclusion.  The proof is complete.

We are ready to present our counterpart of S. Yamada's Theorem 2, \cite{Yam}.

\begin{theorem} \label{onestrata} 
$\Tbar$ is a unique geodesic space.  The
length-minimizing curve connecting points $p, \, q\in\Tbar$ is contained in
the closure of the stratum with label $\Lambda(p)\cap\Lambda(q)$.  The open
segment $\gamma-\{p,\,q\}$ is a solution of the WP geodesic equation on the
stratum with label $\Lambda(p)\cap\Lambda(q)$. For a point $p$ the stratum with
label $\Lambda(p)$ is the union of the length-minimizing open segments
containing $p$.  The closure of each stratum is
a convex set, complete in the induced metric.   
\end{theorem}

\prf  $\Tbar$ is a geodesic space from Proposition \ref{Tcomp}.  For a
length-minimizing curve $\gamma$ we  consider the label behavior of
$\Lambda\circ\gamma$.  From Proposition \ref{finite} $\Lambda\circ\gamma$ only
has a finite number of discontinuities.  From Propositions \ref{refract} and
\ref{corners}, as well as the lifting property of the indicated curves on the
local-manifold covers of $\Mbar$, it follows that $\Lambda\circ\gamma$ is at
most discontinuous at an endpoint.  Since each stratum is a relatively open
subset of its closure in $\Tbar$, it further follows that the value of
$\Lambda$ on the open segment of $\gamma$ precedes its value at each endpoint
(the open segment value is a lower bound for the partial ordering).  It also
follows that $\Lambda(p)\subset\Lambda(q)$ is a necessary condition for $p$ to
be an interior point of a length-minimizing curve  with endpoint $q$.  

A free homotopy class $\alpha$ of a simple closed curve is represented by a
vertex in $\Lambda(p)\cap\Lambda(q)$.   For the  geodesic-length function
$\ell_{\alpha}$, the composition $\ell_{\alpha}\circ\gamma$ is a continuous
function.  The composition vanishes at its domain endpoints and is convex on
its domain interior.  It follows that the composition is identically zero and
consequently that the open segment of $\gamma$ is contained in the stratum
with label $\Lambda(p)\cap\Lambda(q)$.  A stratum is a product of
Teichm\"{u}ller spaces.  The maximal open segment of $\gamma$ is a solution of
the product WP geodesic equation on the stratum.  Consider WP geodesics
$\gamma\,,\gamma'$ parameterized proportional to arc-length by the
unit-interval with common endpoints.  The distance between corresponding
points is a continuous function, vanishing at $0$ and $1$, and convex on
$(0,1)$ from Proposition \ref{discon}.  The distance is identically zero and
the geodesics coincide.  $\Tbar$ is a unique geodesic space.

As note above, $\Lambda(p)\subset\Lambda(q)$ is a necessary condition for  $p$
to be an interior point of a length-minimizing curve with endpoint $q$.  Since
a stratum is a product of Teichm\"{u}ller spaces for which length-minimizing
curves are solutions of the geodesic equation and since solutions can be
extended, it follows that the condition is sufficient for extension.  The
final conclusion follows since $\Tbar$ is a geodesic space.  The closure of a
stratum is convex from the above description of geodesics. The closure of a
stratum is complete in the induced metric from  the completeness of $\Tbar$.
The proof is complete.

We are ready to present the basic result.  $CAT(0)$ is a generalized condition
for a non-positively curved, uniquely geodesic space, \cite[Chap. II.1]{BH}.
With the above result there is little further need to distinguish between
length-minimizing curves and solutions of the geodesic differential equation.
We now also refer to length-minimizing curves parameterized proportional to
arc-length as {\em geodesics}.  

\begin{theorem} \label{catzero} $\Tbar$ is a CAT(0) space.  \end{theorem}

\prf  A length-minimizing curve on $\Tbar$ is approximated by WP geodesics on
$\mathcal{T}$ by choosing sequences of points converging to the endpoints, and
considering the joins parameterized on the original interval.  From
Proposition \ref{discon} the joins converge to the designated geodesics, and
for a pair of geodesics the relative distance functions converge. In particular a limit 
of geodesic triangles satisfying the $CAT(0)$ inequality will also satisfy the inequality, \cite[Chap. II.1]{BH}.
Since the WP metric on $\mathcal T$ has negative curvature, geodesic triangles satisfy the $CAT(0)$ inequality 
\cite[Chap. II.1, Remark 1A.8]{BH}. 
The proof is complete.

The local geometry of geodesics on $\Mbar$ differs from that of $\Tbar$. A
product cuspidal metric is not uniquely geodesic.

\begin{proposition} \label{notcatzero} $\Mbar$ is not locally uniquely
geodesic at the compactification divisor and in particular is not locally a
CAT(0) space.  \end{proposition}

\prf We show that the local manifold cover for a neighborhood of a Riemann
surface having a single node is not uniquely geodesic \cite[Chap. II.1, Prop.
1.4]{BH}.  Introduce hyperbolic metric plumbing coordinates.  For a reference
base point $(s',t')=(s_1',\dots,s_m',t')$ off the $s$-axes, we consider the
curves based at $(s',t')$, disjoint from the $s$-axes, and linking the
$s$-axes. The base point $(s',t')$ lifts to a point of the Teichm\"{u}ller
space $\mathcal T$.  From Corollary \ref{unqgeo} for each possible  value of
the linking number there is a corresponding length-minimizing curve in
$\mathcal T$ (minimizing for curves disjoint from the $s$-axes).  We first
bound the lengths of such linking curves. A general comparison curve is
prescribed by the sequence: a radial line segment in the $t$-coordinate,
followed by an integer number of rotations about a $t$-coordinate circle and
finally a radial line segment returning to the base point.  From Corollary
\ref{wpnormal} the comparison curves can be prescribed with length uniformly
bounded by a multiple of $|t'|$.  It follows for $|s'|,|t'|$ small that the
length-minimizing linking curves are all contained in a small neighborhood of
the origin.

We  consider  length bounds involving the linking number.  For a curve with
linking number $n$ and the minimal absolute value  of the coordinate $t$ on
the curve $|t_0|$, then by Corollary \ref{wpnormal} the length of the curve is
at least a uniform multiple of $|n\,t^3_0|$.  It follows that $|t_0|$  is
bounded by $|t'/n|^{1/3}$; it further follows by considering only the
$t$-radial component of the length that the linking curve has length at least
a uniform multiple of $|t'|-|t'/n|^{1/3}$; the desired bound.  A comparison
curve with linking number one is the $t$-coordinate circle of radius $|t'|$;
its length is bounded by a multiple of $|t'|^3$.  We draw a simple conclusion:
there is a length-minimizing linking curve $\gamma$ of minimal length
(presumably with linking number $\pm 1$).

We bisect $\gamma$. Let $p$ denote $(s',t')$ and $q$ denote the
$\gamma$-midpoint. The length of $\gamma$ is $O(|t'|^3)$; from Corollary
\ref{wpnormal} the length of a curve connecting $p$ to the $s$-axes is at
least a multiple of $|t'|$.  With the length bounds and the fact that $\gamma$
is a solution of the geodesic differential equation it now follows that the
segments of $\gamma$ connecting $p$ to $q$ are length-minimizing for the
neighborhood of the origin.  The neighborhood is not uniquely geodesic.  The
proof is complete.

\section{Applications} \label{appl}

We are interested in understanding the flat subspaces of $\Tbar$.  Our purpose
is to understand the flat geodesic simplices and in particular the flat
geodesic triangles.  Consider a geodesic triangle with distinct vertices $o$,
$p$ and $q$. Parameterize proportional to arc-length the sides $\widehat{op}$
and $\widehat{oq}$ by geodesics $\gamma(t)$ and $\gamma'(t)$, $t\in [0,1]$
with $\gamma(0)=\gamma'(0)=o$.  The distance function
$d(\gamma(t),\gamma'(t))$ is an important measure of the triangle.  We also
require a numerical invariant for noded Riemann surfaces.  Let $\nu(R)$ be the
number of components of $R-\{nodes\}$ that are not thrice-punctured spheres.
The maximal dimension of a flat subspace of the stratum corresponding to $R$
is given by $\nu(R)$.   We use the description of flat subspaces to give a
different proof of a Brock-Farb result: the WP metric is in general not
Gromov-hyperbolic, \cite{BF}.  Recall that a metric space $(M,d)$ is {\em
Gromov-hyperbolic} provided there exists a positive number $\delta$ such that
for each geodesic triangle the $\delta$-neighborhood of a pair of sides
contains the third side, \cite{GH}.

\begin{proposition} \label{flatness} On the augmented Teichm\"{u}ller space
$\Tbar$ of genus $g$, $n$ punctured surfaces \begin{enumerate} \item For
geodesics $\gamma,\,\gamma'$ as above and an interior parameter value,
consider that the values of the distance function $d(\gamma(t),\gamma'(t))$
and its supporting linear function coincide.  The interiors of the geodesics $\gamma,\,\gamma'$
then lie on a submanifold of $\Tbar- \mathcal{T}$ given as the Cartesian
product of geodesics from component Teichm\"{u}ller spaces.  \item  For a
stratum corresponding to a noded Riemann surface $R$, the maximal dimension of
a locally Euclidean isometric submanifold is $\nu(R)$. The maximal value of
$\nu$ is $g-1 + \lfloor \frac{g+n}{2}\rfloor$, which is achieved for an
arrangement with $g$ once-punctured tori and $\lfloor \frac{g+n}{2}\rfloor-1$
four-punctured spheres.  \item For $3g-3+n\ge 3$ the Teichm\"{u}ller space
with the WP metric is not Gromov-hyperbolic.  \end{enumerate}
\end{proposition}

\prf  First, a convex function is necessarily linear if it shares a common
interior value with the supporting linear function.  From Bridson-Haefliger
\cite[Chap. I.1, Defn. 1.10 and Chap. II.3, Prop. 3.1]{BH} provided
$d(\gamma(t),\gamma'(t))$ is linear then the comparison angles formed by the
point triples $(\gamma(t),o,\gamma'(t))$ all coincide.  The {\em flat triangle
lemma} of A. D. Alexandrov can now be applied \cite[Chap. II.2, Prop.
2.9]{BH}.  The convex hull of $o$, $p$ and $q$ in $\Tbar$ is consequently
isometric to the convex hull of a Euclidean triangle with the corresponding
side lengths.  An isometry is prescribed.  There is an associated variation of
geodesics $\beta(t,s)$, parameterized proportional to arc-length,  such that
$\beta(t,0)=\gamma(t),\,\beta(t,1)=\gamma'(t),\,\beta(0,s)=o$ and $\beta(1,s)$
lies on $\widehat{pq}$ with $d(p,\beta(1,s))=sd(p,q)$.  By Theorem
\ref{onestrata} it follows for interior parameter values that $\beta(t,s)$
lies in a single stratum; by Proposition \ref{discon} it follows for interior
parameter values that $\beta(t,s)$ is smooth.  The stratum is a product of
Teichm\"{u}ller spaces. We may apply the techniques of Riemannian geometry,
\cite{CE}.  Since the triangle is flat the contribution to (\ref{secondvar})
from the term $\bigl<R(V,T)T,V\bigr>$ is zero.    For a product of
negatively-curved metrics $\bigl<R(V,T)T,V\bigr>$ vanishes only if the
variation fields $V$ and $T$ everywhere have collinear projections to the
tangent spaces of the factors \cite[Chap. 3, Lemmas 39, 58]{BO'N}.  Since the
projections are collinear the triangle also projects to a geodesic segment in
each component Teichm\"{u}ller space.  The desired first conclusion.

Second, as already indicated for a Riemannian product of negatively-curved
metrics the maximal dimension of a flat subspace is the number of factors. The
dimension of a maximal flat is given by $\nu(R)$ since a punctured Riemann
surface has a positive dimensional Teichm\"{u}ller space provided the surface
is not the thrice-punctured sphere. Once-punctured tori have Euler
characteristic $-1$ and $\dim \mathcal T > 0$; the general surface with $\dim
\mathcal T >0$ has Euler characteristic strictly less than $-1$.  The maximal
statement follows.  The desired second conclusion.

Third, for the Euclidean plane, a positive number $\delta$, and a non
degenerate triangle, a large-scaling provides a triangle $\delta$ with a
$\delta$-neighborhoold of a pair of sides omitting an open segment on the
third side.  A stratum corresponding to a noded surface $R$ with $\nu(R)\ge 2$
contains triangles isometric to $\Delta$.  For such a triangle a triple of
points in $\mathcal T$, one close to each vertex, prescribes a triangle with
measurements close to those of $\Delta$. Independent of the length of an edge,
the geodesic triangle in $\mathcal T$ is uniformly close to the triangle
$\Delta$ with distance estimated only by the distance separating corresponding
vertices.  For a suitable triple, a $\delta$-neighborhood of two joining sides
omits an open segment on the third joining side.  A stratum with $\nu(R)\ge 2$
exists provided $\dim \mathcal T\ge 3$.   The proof is complete.

The maximal simplices in $C(F)$ serve an important role for the geometry of
$\Tbar$.  Since thrice-punctured spheres are conformally rigid, a
$3g-4$-simplex $\sigma$ in $C(F)$ corresponds by $\Lambda^{-1}$ to a unique
marked maximally noded Riemann surface $R_{\sigma}$ in $\Tbar$.  The $Mod$
stabilizer of a maximally noded Riemann surface $R_{\sigma}$ is an extension
of a finite group by a rank $3g-3$ Abelian group, the mapping classes of
products of Dehn twists about the elements of $\sigma$.  The maximally noded
Riemann surfaces serve the role of the {\em maximal rank cusps} for the moduli
space.  Brock studied the finite length WP geodesics from a point of $\mathcal
T$ to the marked noded Riemann surfaces, \cite{Brkwpvs}.  The geodesics from a point can
be extended to include their endpoints in $\Tbar$.  As a consequence of the $CAT(0)$ geometry 
the initial unit tangents for the family of geodesics from a point to a stratum provide for a Lipschitz
map from the stratum to the unit tangent sphere.  Accordingly the image of $\Tbar-\mathcal T$ in
each unit tangent sphere has measure zero and thus the infinite length geodesic rays have tangents
dense in each tangent sphere.   Brock's method for approximating infinite length rays by finite
length rays now provides the following.

\begin{theorem}[\cite{Brkwpvs}] 
The geodesic rays from a point of $\mathcal T$ to the maximally noded 
Riemann surfaces have initial tangents dense in the tangent space.
\end{theorem}

The following is a new consequence of the result.

\begin{corollary}
\label{density}
The geodesics connecting maximally noded Riemann surfaces have tangents dense 
in the tangent bundle of $\mathcal T$.
\end{corollary}

\prf  We consider unit-speed geodesics.  Given a unit-tangent $v$ at a point $p$ of 
$\mathcal T$ and a positive 
number $\epsilon$, we proceed to determine an approximating geodesic.  By the 
above Theorem let $\gamma_-$ be a unit-speed geodesic connecting a point $q$, 
representing a maximally noded Riemann surface to $p$, with the final tangent 
$w$ of $\gamma_-$ within $\epsilon$ of $v$.  For a small positive number $\delta$ 
similarly let $\gamma_+$ be a unit-speed geodesic connecting $p$ to a point $r$, 
representing a maximally noded Riemann surface, with the initial tangent of 
$\gamma_+$ within $\delta$ of $w$.  The geodesics $\gamma_-$ and $\gamma_+$ 
form a vertex at $p$ with angle in the interval $[\pi-\delta,\pi]$.  The three 
points $p$, $q$ and $r$ determine a geodesic triangle $\Delta(p,q,r)$ in 
$\Tbar$ for which there is a comparison triangle $\Delta(\bar p,\bar q,\bar r)$ 
in the Euclidean plane.  By \cite[Chap. II.4, Lem. 4.11]{BH} the vertex angles 
for $\Delta(p,q,r)$ are bounded by the corresponding vertex angles for the Euclidean
triangle $\Delta(\bar p,\bar q,\bar r)$.  In particular the vertex angle at $\bar p$ is 
also in the interval $[\pi-\delta,\pi]$.  It follows for $\Delta(\bar p,\bar q,\bar r)$ 
that the angle at $\bar q$ is at most $\delta$.  For $\delta < \pi/2$ it follows that 
$d(\bar q,\bar r)>d(\bar q, \bar p)$ and there is a point $\bar s$ with 
$d(\bar q,\bar s)=d(\bar q,\bar p)$, $\bar s$ on the geodesic segment 
$\widehat{\bar q\bar r}$.  By trigonometry $d(\bar s,\bar p)=2d(\bar p,\bar q)\sin \delta/2$ 
and thus the comparison point $s$ on the geodesic segment $\widehat{qr}$ satisfies 
$d(s,p)\le 2d(p,q)\sin\delta/2$.  Similarly the midpoints $s_{mid}$ of $\widehat{qs}$ and
$p_{mid}$ of $\widehat{qp}$ satisfy $2d(s_{mid},p_{mid})\le d(s,p)$. To summarize the
considerations, the geodesic $\widehat{qr}$ contains a point $s$ that is within distance
$2d(p,q)\sin \delta/2$ of $p$, and the midpoints $s_{mid}$ of $\widehat{qs}$ and $p_{mid}$ of 
$\widehat{qp}$ are within distance $d(p,q)\sin\delta/2$.

From Proposition \ref{discon} the geodesic segments $\widehat{s_{mid}s}$ and $\widehat{p_{mid}p}$
are sufficiently close in the $C^1$-topology for $\delta$ sufficiently small.  We now
choose $r'$ and $\gamma_+'$ to provide that the tangent at $s'$ is within $\epsilon$ of the final
tangent of $\gamma_-$, which in turn was chosen to be within $\epsilon$ of $v$.  The proof
is complete. 

The following is an immediate consequence of the above result.

\begin{corollary} 
\label{convexhull}
$\Tbar$ is the closed convex hull of the subset of marked maximally noded 
Riemann surfaces.  
\end{corollary}

We combine the above and follow the outline of the Masur-Wolf approach to give
an immediate proof of the Masur-Wolf theorem, \cite[Theorem A]{MW}.  The
classification of simplicial automorphisms of the curve complex $C(F)$ by N.
Ivanov \cite{Ivaut}, M. Korkmaz \cite{Krk}, and F. Luo \cite{Luaut} is an
essential consideration.

\begin{theorem}
\label{WPisom} 
For $3g-3n>1$ and $(g,n)\ne (1,2)$, every WP isometry of
$\mathcal T$ is induced by an element of the extended mapping class group.
\end{theorem}

\prf  An isometry of $\mathcal T$ extends to an isometry of the completion
$\Tbar$.  By Theorem \ref{onestrata} an isometry of $\Tbar$ necessarily
preserves the strata structure and the incidence relations.  It follows that
an isometry induces a simplicial automorphism of $C(F)$.  From the results of
Ivanov \cite{Ivaut}, Korkmaz \cite{Krk} and Luo \cite{Luaut}, every simplicial
automorphism coincides with the induced automorphism of an extended mapping
class.  In particular for a WP isometry of $\Tbar$, there is an extended
mapping class such that the two mappings coincide on the subset of maximally
noded Riemann surfaces.  The conclusion now follows from Corollary
\ref{convexhull}.  The proof is complete.

A complete, convex subset $\mathcal C$ of a $CAT(0)$ space is the base for an
{\em orthogonal projection}, \cite[Chap. II.2, Prop. 2.4]{BH}.  A fibre of the
projection is the unique geodesic realizing the distance between its points
and the base.  The projection is a retraction that does not increase distance.
The distance $d_{\mathcal C}$ to $\mathcal C$ is a convex function satisfying
$|d_{\mathcal C}(p)-d_{\mathcal C}(q)|\le d(p,q)$, \cite[Chap. II.2, Prop.
2.5]{BH}.  Examples of complete, convex sets $\mathcal C$ are: points,
complete geodesics, and fixed-point sets of isometry groups.  Furthermore in
the case of $\Tbar$ with Theorem \ref{onestrata} the closure of each
individual stratum is the base of a projection (local projections are
prescribed on $\Mbar$).  On $\Tbar$ a tubular neighborhood of an   stratum is
fibered by the projection-geodesics.  For the local understanding of the
distance to the stratum  we now consider a refinement of the prescription for
a product cuspidal metric. In particular by Corollary \ref{wpnormal} the WP
metric has an expansion with an order-three approximation about the $x$-axes   

$$ ds^2=(d\mu^2+dr^2+r^6d\theta^2)(1+O(\|r\|^3). $$

\begin{corollary} 
\label{distD}
For a stratum $\sigma$ defined by vanishing of
the geodesic-length sum $\ell=\ell_1+\cdots+\ell_n$, the distance to the
stratum is given locally as  $d(p,\sigma)=(2\pi\,\ell)^{1/2}+O(\ell^2)$.
\end{corollary}

\prf We begin by considering distances on $\Mbar$.  For a prescribed stratum
and point choose a local-manifold cover $\widehat{\mathcal{U}}$ with the point
corresponding to the origin and (image) stratum $\sigma$ corresponding to the
$x$-axes for the normal form of $ds^2$.  The distance to the stratum
$d(p,\sigma)$ is estimated from above by considering the radial line from a
point $p$ to the origin.  The bound is $\|p_y\|+O(\|p_y\|^4)$ in terms of the
$y$-projection of $p$.  We next consider a lower bound for $d(p,\sigma)$.  A
curve in $\widehat{\mathcal {U}}$ connecting $p$ to $\sigma$ can be
isometrically lifted to $\Tbar$; it suffices to examine curves that are
length-minimizing on $\Tbar$.  From \cite[Chap. II.2, Prop. 2.4]{BH} for $p$
close to the origin there is a geodesic $\gamma\subset\widehat{\mathcal{U}}$
connecting $p$ to $q\in \sigma$ and $\gamma$ provides the length-minimizing
curve connecting to $\sigma$ for each point of $\gamma$.  For $p'\in \gamma$,
as noted, $d(p',\sigma)$ is bounded above by the Euclidean norm of $p_y'$.
Since $ds^2$ is likewise bounded below in terms of $dr^2$ it follows that
$d(p',\sigma)$ is actually comparable to $\|p'_y\|$ and consequently that the
maximum Euclidean $y$-height of the $\gamma$-segment $\widehat{p'q}$ is
likewise bounded.  It now follows overall that $$
d(p',\sigma)=\int_{\widehat{p'q}}(d\mu^2+dr^2+r^6d\theta^2)^{1/2}+O(\|p'\|^4).
$$ The explicit integral is minorized by choosing the radial line in the
$y$-coordinate; the resulting lower bound is $\|p'_y\|$.  Thus the distance to
the stratum on $\Mbar$ and on $\Tbar$ is $d(p,\sigma)=\|p_y\|+O(\|p_y\|^4)$.
In \cite[Example 4.3]{Wlhyp} a relationship is provided for the
geodesic-length functions and the present hyperbolic metric plumbing
coordinates; it is shown that $$ \ell_j=\frac{2\pi^2}{-\log |t_j|}+
O\bigr(\frac{1}{(\log^2|t_j|)}\sum\limits_{k=1}^n\frac{1}{(\log^2|t_k|)}
\bigr).  $$ We may rearrange terms and substitute the relation
$\varrho_j^2=\frac{4\pi^3}{-\log |t_j|}$ to obtain the expansion $$
\varrho_j^2=2\pi\ell_j(1+O\bigr(\ell_j\sum\limits_k\ell_k^2)).  $$ The
distance expansion now follows from Corollary \ref{wpnormal}. The proof is
complete.

The WP gradients of the geodesic-length functions also have general
expansions, \cite[II: Sec. 2.2, Lemmas 2.3 and 2.4]{Wlspeclim}.  We have for
$\epsilon$ positive and geodesic-length functions
$\ell_{\alpha}\,\le\,\ell_{\beta}\,\le\,\epsilon$ that 

\begin{gather*} \bigl<\grad \ell_{\alpha},\grad \ell_{\alpha}\bigr> =
\frac{2}{\pi}\ell_{\alpha}\ +\ O(\ell^3_{\alpha}) \\ \bigl<\grad
\ell_{\alpha},\grad\ell_{\beta}\bigr> = O(\ell_{\alpha}\ell_{\beta}^2)
\end{gather*} 
with  constants independent of the surface, or in particular for
$\lambda_*=(2\pi\ell_*)^{1/2}$ that 

\begin{gather*} \bigl<\grad \lambda_{\alpha},\grad\lambda_{\alpha} \bigr>=1\
+\ O(\lambda_{\alpha}^4)\\ \bigl<\grad \lambda_{\alpha},\grad\lambda_{\beta}
\bigr>= O(\lambda_{\alpha}\lambda_{\beta}^3).  \end{gather*}

We consider applications of the gradient bounds.  We now consider simplices
$\sigma$ and $\tau$ with free homotopy classes $\alpha$ in $\sigma$ and
$\beta$ in $\tau$ with (all) representatives intersecting.  The simplice
$\sigma$ and $\tau$ do not precede a common simplex.  There is also a positive
lower bound for the corresponding geodesic-length sum
$\ell_{\alpha}+\ell_{\beta}$.   In particular from the {\em collar result}
there is a positive constant $\ell_0<2$ such that about a geodesic $\alpha$
with $\ell_{\alpha}\le\ell_0$, there is an embedded collar of width $2\log
2/\ell_{\alpha}$: in which case $\ell_{\beta}$ is at least the width,
\cite{Rn}.    Consider a WP length-minimizing curve $\gamma$ connecting the
strata $\mathcal{S}(\sigma)$ and $\mathcal{S}(\tau)$.  On the curve $\gamma$
the geodesic-length functions  $\ell_{\alpha}$ and $\ell_{\beta}$ each vanish
($\ell_{\alpha}$  on $\mathcal S(\sigma)$; $\ell_{\beta}$  on $\mathcal
S(\tau)$) and are each unbounded.    The following is now the consequence of
the universal bounds for the gradients. 

\begin{corollary} \label{posdist} There is a positive constant $\delta$ such
that null strata $\mathcal{S}(\sigma)$ and $\mathcal{S}(\tau)$ either have
intersecting closures or $d(\mathcal{S}(\sigma),\mathcal{S}(\tau))\ge \delta$.
\end{corollary}

We also recognize from the relations  for a stratum $\sigma$, defined by the
vanishing of the geodesic-length sum $\ell=\ell_1+\cdots+\ell_n$, that the
vector fields $\{\grad\lambda_1,\dots,\grad\lambda_n\}$ play the role of {\em
normal Fermi fields} \cite[Sec. 2.3]{Gray}.  Recall in particular for a
tubular neighborhood of a submanifold in a Riemannian manifold, Fermi
coordinates are the relative analog of {\em normal} coordinates about a point.
First consider a neighborhood $\mathcal N$ of the $0$-section of the normal
subbundle of the tangent bundle of the submanifold.   With the restriction of
the  exponential map $\mathcal N$ is identified with a tubular neighborhood of
the submanifold.  An orthonormal frame for the normal bundle provides Fermi
coordinates for the fibres of $\mathcal N$, and also for a tubular
neighborhood of the submanifold upon composition with the inverse of the
exponential map.  The unit-speed geodesics normal to the submanifold are given
in terms of the Fermi coordinates as the unit-speed linear rays from the
$0$-section.  The tangent fields of the unit-speed geodesics normal to the
submanifold are constant sums of the Fermi coordinate tangent fields (the {\em
normal Fermi fields}).

We are ready to present the analogy.  For constant sums of the vector fields
$\grad\lambda_*$,  the WP distance between endpoints of integral curves nearly
equals the integral curve length.    Consider for a positive vector
$c=(c_1,\dots,c_n)$ the integral curves of
$\mathfrak{v}=-\sum_j\,c_j\grad\lambda_j$.  The time-one integral curve of
$\mathfrak{v}$ with initial point
$2\pi(\ell_1,\dots,\ell_n)=(c_1^2,\dots,c_n^2)$ has terminal point at distance
$O(\|c\|^4)$ to $\sigma$. Further from Corollary \ref{distD} the point
$2\pi(\ell_1,\dots,\ell_n)=(c_1^2,\dots,c_n^2)$ is at WP distance
$\|c\|+O(\|c\|^4)$ from $\sigma$.  From the  gradient relations above we have
that $\|\mathfrak{v}\|_{WP}=\|c\|+O(\|c\|^4)$ and that the time-one integral
curve has the same WP length.   For $\|c\|$ small, the time-one integral
curves of $\mathfrak v$ have endpoints at distance nearly equal to the curve
length.  The integral curves   approximate WP geodesics.   The integral curves
of $\mathfrak{n}=\sum_j\ell^{-1/2}\lambda_j\grad \lambda_j$ also have length
nearly equal to the distance between endpoints.   The time $(2\pi\ell)^{1/2}$
integral curve of $\mathfrak n$ connects $(\ell_1,\dots,\ell_n)$ and $\sigma$;
for $\ell$ small $\mathfrak{n}$ is approximately the WP unit normal field to
$\sigma$. Also by Corollary \ref{distD} $\mathfrak n$ approximates $\grad
d_{WP}(\cdot,\sigma)$.

\section{The structure of geodesic limits} \label{geolim}

We investigate  sequences of geodesics.  $\Tbar$ is a complete metric space
with a compact quotient $\Mbar$.  We anticipate that  the compactness is
manifested in the structure of the space of geodesics for $\Tbar$.  We find
that geometric limits of geodesics are described by {\em polygonal paths} and
products of Dehn twists.  Specifically for a sequence of bounded length
geodesics there is a subsequence of $Mod$-translates that converges
geometrically (sequences of products of Dehn twists are applied to
subsegments) to a  polygonal path, a piecewise geodesic curve connecting
strata. We consider an application of the result and show that each
fixed-point free element of the mapping class group $Mod$  has a geodesic axis
in $\Tbar$; the axis is unique and lies in $\mathcal T$  when the element is
irreducible.  Furthermore  irreducible elements have either coinciding or
divergent axes.  The present results provide a different approach for the
considerations of G. Daskalopoulos and R. Wentworth \cite{DW2}.

Sequences of geodesics can have special behavior for product cuspidal metrics.
We present an example.  Consider the half-plane $\mathbb{R}_{\ge
0}\times\mathbb{R}$ with coordinates $(r,\theta)$ and the identification space
$\mathbb{R}_{\ge 0}\times\mathbb{R}/\{(0,\theta)\sim(0,\theta')\}$ with {\em
metric} $dr^2+r^6d\theta^2$.  Denote the special point
$\{(0,\theta)\sim(0,\theta')\}$ by $O$.  For the isometry
$T:(r,\theta)\rightarrow(r,\theta+1)$ consider the unit-speed geodesics
$\gamma_n$ connecting $(r_0,\theta_0)$ and $T^n(r_1,\theta_1)$.  For $n$ large
the length of $\gamma_n$ is nearly $r_0+r_1$; we can provide that the
$\gamma_n$ are essentially parameterized on the interval $[0,r_0+r_1]$.  By
elementary considerations of differential equations, on the parameter interval
$[0,r_0]$ the sequence $\{\gamma_n\}$ converges to the $\theta=\theta_0$ line
segment $\widehat{(r_0,\theta_0)O}$.  On the parameter interval
$[r_0,r_0+r_1]$ the sequence $\{T^{-n}\gamma_n\}$ converges to the
$\theta=\theta_1$ line segment $\widehat{O(r_1,\theta_1)}$.  In effect the
geodesic sequence $\{\gamma_n\}$ is described by the polygonal path
$\widehat{(r_0,\theta_0)O}\cup\widehat{O(r_1,\theta_1)}$ and the sequence of
transformations $\{T^n\}$.  Furthermore the  curve
$\widehat{(r_0,\theta_0)O}\cup T^n\widehat{O(r_1,\theta_1)}$ is continuous and
has distance in the sense of parameterized unit-speed curves to $\gamma_n$
that tends to zero as $n$ tends to infinity.

We consider the local description of the map from $\Tbar$ to $\Mbar$.
Associated to a $k$-simplex $\sigma$ is the rank $k+1$ Abelian group
$Mod_{\sigma}$ of mapping classes of products of Dehn twists about the
elements of $\sigma$.  $Mod_{\sigma}$ stabilizes the $\sigma$-null stratum
$\mathcal{S}(\sigma)$.  For a point $p\in \mathcal{S}(\sigma)$ the stabilizer
$Mod(p)\subset Mod$ is a group extension of a finite group $G(p)$ by
$Mod_{\sigma}$.  Furthermore for a point $p$ we can prescribe a suitable basis
$\{\mathcal{U}\}$ of $Mod(p)$ invariant neighborhoods. Neighborhoods of the
projection of $p$ to $\Mbar$ are given as $\mathcal{U}/Mod(p)$.  Furthermore
each quotient $\mathcal{U}/Mod_{\sigma}$ is a local manifold cover.  We can
further prescribe that the quotients $\mathcal{U}/Mod_{\sigma}$ are relatively
compact in a fixed quotient $\mathcal{U}'/Mod_{\sigma}$; the quotients
$\partial\mathcal{U}/Mod_{\sigma}$ are accordingly compact.

Since $\Mbar$ is compact, given a sequence of points of $\Tbar$  there exists
a subsequence and associated elements of $Mod$ such that the sequence of image
points converges.  Accordingly we consider sequences of unit-speed
parameterized geodesics with initial points converging.

\begin{proposition} \label{geolims} Consider a sequence of unit-speed
geodesics $\{\gamma_n'\}$ with initial points converging to $p_0$, lengths
converging to a positive value $L'$ and parameter intervals converging to
$[t',t'']$ with $L'=t''-t'$.  There exists an associated partition
$t'=t_0<t_1<\cdots <t_k=t''$ of the interval; simplices
$\sigma_0=\Lambda(p_0),\sigma_1,\dots,\sigma_k$; and points
$p_1\in\mathcal{S}(\sigma_1),\dots ,p_k\in\mathcal{S}(\sigma_k)$ on the null
strata.  

The data satisfies $L(\widehat{p_jp_{j+1}})=t_{j+1}-t_j$ for
$j=0,\dots,k-1$ and for the stratum with label
$\tau_j=\Lambda(p_j)\cap\Lambda(p_{j+1})$  then:
$\tau_0$ strictly precedes $\sigma_1$ if $k>1$; $\tau_{k-1}$ strictly precedes
$\sigma_{k-1}$ if $k>1$; $\tau_j$ strictly precedes $\sigma_j$ and
$\sigma_{j+1}$ for $j=1,\dots,k-2$.  The concatenation of geodesic segments
$\widehat{p_0p_1}\cup\widehat{p_1p_2}\cup\dots\cup\widehat{p_{k-1}p_k}$ is the
unique length-minimizing curve connecting $p_0$ to $p_k$ and intersecting in
order the closures of the strata
$\mathcal{S}(\sigma_1),\dots,\mathcal{S}(\sigma_{k-1})$.

There is a subsequence $\{\gamma_n\}$ of the geodesics and sequences of products of Dehn
twists $T_{(j,n)}\in Mod_{\sigma_j-\tau_j}$, $j=0,\dots,k-1$, such that on the
parameter interval $[t_j,t_{j+1}]$ the  geodesic segments $T_{(j,n)}\circ\cdots\circ
T_{(0,n)}\gamma_n$ converge to $\widehat{p_jp_{j+1}}$ in the sense of
parameterized unit-speed curves.  Furthermore the distance between the
parameterized unit-speed curves $\gamma_n$ and
$\widehat{p_0p_{(k,n)}},\,p_{(k,n)}=(T_{(k-1,n)}\circ\cdots\circ
T_{(0,n)})^{-1}p_k$, tends to zero for $n$ tending to infinity.  The sequence
of transformations $\{T_{(0,n)}\}$ is either trivial or unbounded.  The
sequences of transformations $\{T_{(j,n)}\}$, $j=1,\dots,k-1$ are unbounded.
\end{proposition}

\prf  The main argument is to provide the two steps for determining the
individual geodesic segments $\widehat{p_jp_{j+1}}$.  The overall argument is then a
finite induction.  For the first step choose a neighborhood $\mathcal{U}$ of
$p_0$ with $\partial\mathcal{U}/Mod_{\sigma_0}$ compact.  For each geodesic
$\gamma_n'$ let $q_n$ be the first point of intersection with
$\partial\mathcal{U}$.  Either a subsequence of the points $q_n$ converges to
a point $q'$ or we select elements $T_{(0,n)}\in Mod_{\sigma_0}$ such that the
images $T_{(0,n)}q_n$ lie in a relatively compact fundamental domain for the
action of $Mod_{\sigma_0}$  on $\partial\mathcal{U}$.  For the situation of
selecting elements $T_{(0,n)}$ there is a subsequence $\{T_{(0,n)}q_n\}$
convergent to a point $q'$ and the sequence $\{T_{(0,n)}\}$ is unbounded.  Now
the group $Mod_{\sigma_0}$ fixes $p_0$ and a sequence of points converges to
$q'$.  The group $Mod_{\sigma_0}$ is a direct product with a factor
$Mod_{\sigma_0-\tau_0}$ for $\tau_0=\Lambda(q')$; for $\{T_{(0,n)}q_n\}$
converging to $q'$ it is a basic feature of the $\Tbar$ topology that the
$T_{(0,n)}$ can be replaced with their $Mod_{\sigma_0-\tau_0}$ factors and the
resulting sequence also converges.  Finally since geodesics in a $CAT(0)$
space depend continuously on endpoints \cite[Chap. II.1, Prop. 1.4]{BH},  the
appropriate geodesic segments from $p_0$ to $q_n$ or to $T_{(0,n)}q_n$
converge to $\widehat{pq'}$, as claimed.

The second step is to show for a subsequence of the geodesics that maximal
initial segments converge to segments of the prolongation of $\widehat{pq'}$
in the stratum $\tau_0$.  We now write $T_{(0,n)}\gamma_n'$ whether the
transformations are trivial or not.  We preliminarily note that the
subsequence necessarily converges on a closed interval.  In particular for a
subsequence converging uniformly to a segment $\widehat{p_0q''}-\{q''\}$ and
$\epsilon$ small, consider the point on $\widehat{p_0q''}$ distance $\epsilon$
before $q''$; a corresponding sequence of points, one on each
$T_{(0,n)}\gamma_n'$, is determined.  On each $T_{(0,n)}\gamma'_n$ consider
the point $\epsilon$ further along than the referenced point; the resulting
points are at distance $3\epsilon$ from $q''$ for $n$ large.  The interval of
convergence is indeed closed.  There are now three possibilities for the
interval: i) a subsequence $\{T_{(0,n)}\gamma'_n\}$ converges on the entire parameter
interval $[t',t'']$ to a segment of the forward prolongation of
$\widehat{p_0q'}$ (and the overall convergence argument is complete), ii)  a
subsequence $\{T_{(0,n)}\gamma'_n\}$ converges on  $[t',t_1]$, $t_1<t''$, to
$\widehat{p_0p_1}-\{p_1\}$ and $\Lambda(p_1)$ properly succeeds $\tau_0$, or
iii)  a subsequence $\{T_{(0,n)}\gamma_n'\}$ converges on $[t',t_*]$,
$t_*<t''$, to $\widehat{p_0q''}$, with $\widehat{p_0q''}$ having a nontrivial
forward prolongation in the stratum $\tau$.     We  examine case iii).  We
examine the behavior of the subsequence $\{T_{(0,n)}\gamma_n'\}$ in a
neighborhood of $q''$.  We can again apply the argument from the beginning of
the proof to determine elements $S_n\in Mod_{\tau_0}$ such that a subsequence
of $\{S_n\circ T_{(0,n)}\gamma'_n\}$ converges in a neighborhood of $q''$.
The limit is length-minimizing.  From Proposition \ref{refract} the
subsequence $\{S_n\circ T_{(0,n)}\gamma_n'\}$ converges to the prolongation of
$\widehat{p_0q''}$; from the above observation concerning the $\Tbar$ topology
and the $Mod_{\tau_0}$ action the subsequence $\{T_{(0,n)}\gamma_n'\}$
also converges to the prolongation.  We now summarize the convergence
considerations: for a maximal parameter interval of convergence of a
subsequence $\{T_{(0,n)}\gamma'_n\}$ either: i) the parameter interval is
$[t',t'']$, or ii) the interval is $[t',t_1]$, $t_1<t''$, and the limit is
$\widehat{p_0p_1}$ with $\Lambda(p_1)$ strictly succeeding
$\Lambda(p_0)\cap\Lambda(p_1)$.    

We now proceed and apply the considerations of the above two paragraphs to the
subsequence $\{T_{(0,n)}\gamma_n'\}$ considered on the interval $[t_1,t'']$.
The initial points converge to $p_1$.  Elements $T_{(1,n)}\in
Mod_{\sigma_1-\tau_1}$ are determined such that a further subsequence
$\{T_{(1,n)}\circ T_{(0,n)}\gamma_n'\}$ converges as in case i) or case ii).
A stratum $\tau_1$ is prescribed.  The simplex $\sigma_1$ properly succeeds
$\tau_1$ for otherwise (from the observation concerning the $\Tbar$ topology
and the $Mod_*$ action) elements $T_{(1,n)}$ are not required and the
subsequence $\{T_{(0,n)}\gamma_n'\}$ on $[t',t_1+\epsilon)$ converges to a
curve that is not length-minimizing by Proposition \ref{refract}.  We now note
that $Mod$ preserves the strata structure of $\Tbar$.  Since the entire
considerations including the initial-point convergence can be applied to the
sequence $\{\gamma_n'\}$ {\em starting} from an arbitrary value $t'''$ and
proceeding in the negative $t$-direction, we  observe that consequently a
finite partition $t'=t_0<t_1<\dots<t_k=t''$ is determined.   Points
$p_0,\dots,p_k$; strata $\sigma_0,\dots,\sigma_{k}, \tau_0,\dots,\tau_{k-1}$;
and sequences $\{\gamma_n\},\,\{T_{(0,n)}\},\dots,\{T_{(k-1,n)}\}$ are
determined.  The desired properties are provided in the above construction
with only two remaining matters: the sequences $\{T_{(j,n)}\}$,
$j=1,\dots,k-1$, are unbounded and the length-minimizing property of the
concatenation of the points $p_0,\dots,p_k$.

We consider the length property first.  A candidate length-minimizing curve is
given as a concatenation
$\mathcal{C}=\widehat{p_0q_1}\cup\widehat{q_1q_2}\cup\dots
\cup\widehat{q_{k-1}p_k}$  with $p_0\in \mathcal{S}(\sigma_0)$ and
$q_j\in\mathcal{S}(\sigma_j)$ for $j=1,\dots,k-1$.  Since the group
$Mod_{\sigma_j-\tau_j}$ stabilizes $\overline{\mathcal{S}(\sigma_j)}$ it
follows that the concatenation
$\mathcal{C}_T=T^{-1}_{(0,n)}\widehat{p_0q_1}\cup T^{-1}_{(0,n)}\circ
T^{-1}_{(1,n)}\widehat{q_1q_2}\cup\cdots\cup T^{-1}_{(0,n)}\circ\cdots\circ
T^{-1}_{(k-1,n)}\widehat{q_{k-1}p_k}$ is a continuous curve connecting $p_0$
to $p_{(k,n)}=(T_{(k-1,n)}\circ\cdots\circ T_{(0,n)})^{-1}p_k$.   From the
overall construction $d(\gamma_n(t''),p_{(k,n)})$ tends to zero; from
\cite[Chap. II.1, Prop. 1.4]{BH} the distance between $\gamma_n$ and
$\widehat{p_0p_{(k,n)}}$ is consequently small for $n$ large, as claimed.  The
three curves $\gamma_n$, $\mathcal{C}_T$ and $\widehat{p_0p_{(k,n)}}$ each
approximately connect $p_0$ and $p_{(k,n)}$.  Since
$L(\mathcal{C})=L(\mathcal{C}_T)$, it follows that $L(\mathcal{C})\ge
\lim\limits_n L(\gamma_n)=L'$. Thus each suitable concatenation has length at
least $L'$ with the minimum achieved for the arrangement of points
$\{p_0,p_1,\dots,p_k\}$.  It remains for $k>1$ to establish uniqueness.
Consider geodesics $\alpha_j(s)$, parameterized  on the unit-interval,
connecting $p_j$ to $q_j$, $j=1,\dots,k-1$.  The concatenation
$\widehat{p_0\alpha_1(s)}\cup\widehat{\alpha_1(s)\alpha_2(s)}\cup\dots\cup\widehat{\alpha_{k-1}(s)p_k}$
satisfies the strata hypothesis and by Theorem \ref{catzero} has length a
convex function of the parameter $s$.  Now for a $CAT(0)$ space either the
distance from a point to a geodesic is a strictly convex function, or the
point lies on the prolongation of the geodesic, \cite[Chap. II.1, Defn.
1.1]{BH}.  The points $\{p_0,p_1,q_1\}$ do not lie on a common geodesic since
$\tau_0=\Lambda(p_0)\cap\Lambda(p_1)$ strictly precedes $\sigma_1$ and
$\sigma_1\subset\Lambda(p_1)\cap\Lambda(q_1)$.  In consequence either
$L(\widehat{p_0\alpha_1(s)})$ is a strictly convex function or $\alpha_1(s)$
is constant.  It follows that the length of a non constant family is a
strictly convex function.  The minimal-length concatenation is  unique.  

The final matter is the unboundedness of the sequences $\{T_{(j,n)}\}$.  The
elements $T_{(j,n)}$ are selected to provide that the initial curve segments
lie in relatively compact sets.  We can prescribe that each sequence
$\{T_{(j,n)}\}$ is either trivial or unbounded.  In a neighborhood of the
parameter value $t_j$ the concatenation
$\mathcal{C}_j=\widehat{p_{j-1}p_j}\cup T^{-1}_{(j,n)}\widehat{p_jp_{j+1}}$ is
the approximation to $T_{(j-1,n)}\circ\cdots\circ T_{(0,n)}\gamma_n$.  Since
$\gamma_n$ is length-minimizing, it follows for $n$ large that the
concatenation $\mathcal{C}_j$ is arbitrarily close to length-minimizing.  If
$\{T_{(j,n)}\}$ were trivial then $\widehat{p_{j-1}p_j}\cup
\widehat{p_jp_{j+1}}$ would be length-minimizing in contradiction of Theorem
\ref{onestrata}, since $p_j\in\overline{\mathcal{S}(\sigma_j)}$ and $\sigma_j$
strictly succeeds the stratum of one of the connecting geodesic segments.  The proof is
complete.

We now introduce terminology for the data for a convergent sequence of
geodesics.  Consider as above a convergent sequence $\{\gamma_n\}$ with data
$\{p_j\}$ and $\{T_{(j,n)}\}$.  

\begin{definition} \label{polygon} For a convergent sequence of geodesics
$\{\gamma_n\}$ the vertices are the points $p_j$, $j=0,\dots,k$; the vertex
concatenation is
$\widehat{p_0p_1}\cup\widehat{p_1p_2}\cup\cdots\cup\widehat{p_{k-1}p_k}$ and
an approximating concatenation is $T^{-1}_{(0,n)}\widehat{p_0p_1}\cup
T_{(0,n)}^{-1}\circ T_{(1,n)}^{-1}\widehat{p_1p_2}\cup\cdots\cup
T_{(0,n)}^{-1}\circ\cdots\circ T_{(k-1,n)}^{-1}\widehat{p_{k-1}p_k}$.
\end{definition} 

We are ready to consider the matter of existence of axes for elements of
$Mod$.  We present a different approach towards certain results of G.
Daskalopoulos and R. Wentworth \cite{DW2}.  They show that each irreducible
(pseudo Anosov) mapping class has a unique invariant axis and that non
commuting irreducible mapping classes have divergent axes.  To provide a
context we first recall the Thurston-Nielsen classification of mapping classes
\cite[Expos\'{e}s 9, 11]{TrvT}.  A mapping class is {\em irreducible} provided
no power fixes the free homotopy class of a simple closed curve.  A mapping
class is precisely one of: periodic, irreducible or reducible, \cite{TrvT}.
Reducible classes are first analyzed in terms of mappings of proper
subsurfaces.  For a reducible mapping class $[h]$ an invariant is
$\sigma_{[h]}$ the maximal simplex fixed by a power.  A general invariant of a
transformation $S$ is its translation length: $\inf\limits_p d(p,Sp)$.

\begin{theorem} \label{axes} A mapping class S acting on $\Tbar$ either has
fixed-points or positive translation length realized  on a closed, convex set
$\mathcal{A}_S$, isometric to a metric space product $\mathbb{R}\times Y$.
In the latter case the isometry $S$ acts on  $\mathbb{R}\times Y$ as the
product of a translation of $\mathbb{R}$ and $id_Y$.  For $S$ irreducible the
translation length is positive and $\mathcal{A}_S$ is a geodesic in
$\mathcal{T}$.  For $S$ reducible the null stratum $\mathcal{S}(\sigma_{[h]})$
is a product of Teichm\"{u}ller spaces
$\prod\mathcal{T}'\times\prod\mathcal{T}''$ with a power $S^m$ fixing the
factors, acting by a product of: irreducible elements $S'$ on $\mathcal{T}'$
with axis $\gamma_{S'}$ and the identity on each $\mathcal{T}''$;
$\mathcal{A}_S\subset\prod\gamma_{S'}\times\prod\Tbar''$.  \end{theorem}

\prf  The first matter is to establish that the translation length is
realized.  We consider a sequence of geodesics $\{\gamma_n'\}$ parameterized
on $[a,b]$ with $S$ connecting endpoints, $S(\gamma_n'(a))=\gamma_n'(b)$, and
$\lim\limits_n L(\gamma_n')$ the translation length.  Apply elements of $Mod$
and according to Proposition \ref{geolims} select a convergent subsequence
$\{\gamma_n\}$ with vertices $\{p_0,\dots,p_k\}$; each $\gamma_n$ has
endpoints connected by a conjugate of $S$. For the special situation of
translation length zero then $a=b$ and the vertex concatenation is the
singleton $\{p_0\}$.  The main matter is to determine the distance between the
$Mod$ orbit of $p_0$ and that of $p_k$.  For an $\epsilon$ approximating
concatenation to $\gamma_n$ and a conjugate $Q$ we have
$Q(\gamma_n(a))=\gamma_n(b)$ with $\gamma_n(a)$ within $\epsilon$ of $p_0$ and
$\gamma_n(b)$ within $\epsilon$ of the prescribed endpoint $p_{(k,n)}$; it
follows that $d(Q(p_0),p_{(k,n)})<2\epsilon$.  It follows that the distance
between the $Mod$ orbit of $p_0$ and that of $p_k$ is zero.  The distance is
also given by considering $Mod$ translates of $p_k$ in a neighborhood of
$p_0$.  From the preliminary discussion there is a positive lower bound for
the distance between the points of the $p_k$ orbit.  It now follows that for a
suitable $\epsilon$ and $n'$ above, the distance inequality implies that
$Q(p_0)=p_{(k,n')}$.  For the value $n'$ the approximating concatenation
connects $p_0$ and $p_{(k,n')}$ and has length $\lim\limits_n L(\gamma_n)$,
the minimal translation length.  It  follows that there is only one geodesic segment in
the concatenation and  that $Q$ realizes its translation length for the
segment $\widehat{p_0Q(p_0)}$. $S$  realizes its  translation length for an
image of the segment.

By general considerations for positive translation length $S$ realizes its
translation length on axes in $\Tbar$, \cite[Chap. II.6, Defn 6.3, Thrm.
6.8]{BH}.  An axis is isometric to $\mathbb{R}$ and may not be unique; we
consider specifics.  $S$ stabilizes each axis and thus stabilizes the stratum
of an axis.  Since in fact an irreducible mapping class only stabilizes the
single stratum $\mathcal{T}$, it follows that an irreducible class is
fixed-point free with axes in $\mathcal{T}$.  Since axes are parallel and the
distance between geodesics in a Teichm\"{u}ller space is a strictly convex
function, it  follows that an irreducible axis is unique in $\Tbar$.  Now for
a reducible mapping class $S$ a power $P$ fixes a simplex if and only if the
power fixes the corresponding null stratum.  For the maximal simplex
$\sigma_S$ the geodesic-length sum
$\mathcal{L}_S=\sum_{\alpha\in\sigma_S}\ell_{\alpha}$ restricts on each
geodesic of $\Tbar$ to either the zero-function or a strictly convex function.
Since $\mathcal{L}_S$ is $P$-invariant on an axis of $P$ the restriction is
the zero-function.  Thus the axes of $P$ are contained in the null stratum
$\mathcal{S}(\sigma_S)$.  On considering a power of $P$ we can further arrange
that $P^m$ stabilizes the components $F-\cup_{\alpha\in\sigma_S}\{\alpha\}$
of the reference surface, and by \cite[Expos\'{e} 11, Thrm 4.2]{TrvT} that on
each component of $F-\cup_{\alpha\in\sigma_S}\{\alpha\}$ the restriction of
$P^m$ is either the trivial  or an irreducible mapping class.  For positive
translation length at least one factor is an irreducible mapping class since
by \cite[Chap. II.6, Thrm. 6.8(2)]{BH} the translation length of the power is
also positive. It follows that $P^m$ realizes its translation length on a
product $\mathcal{A}_{P^m}=\prod\gamma'\times\prod\Tbar''$ contained in the
closure of $\mathcal{S}(\sigma_S)$.  By \cite[Chap. II.6, Thrm. 6.8(4)]{BH}
$\mathcal{A}_{P^m}$ is isometric to a product $\mathbb{R}\times Y$ with $S$
stabilizing the product and acting thereon by the product of a translation and
a periodic element.  $\mathcal{A}_{P^m}$ is itself a complete $CAT(0)$ space.
It follows that  $Y$ is a complete $CAT(0)$ space and by \cite[Chap. II.2,
Cor. 2.8]{BH} that a periodic element has a non-empty closed convex
fixed-point set.  The proof is complete.

We recall  notions regarding the behavior of geodesic rays.  Unit-speed  rays
$\gamma(t)$, $\gamma'(t)$ are {\em asymptotic} provided the limit
$\lim\limits_{t\rightarrow\infty}d(\gamma(t),\gamma'(t))$ is zero and are {\em
divergent} provided the limit is infinity.   For values $t$, $\tau$ the points
$\gamma(0)$, $\gamma'(0)$, $\gamma'(\tau)$ and $\gamma(t)$ determine a
quadrilateral; on comparing side lengths we have that $|\tau-t|\le
d(\gamma(0),\gamma'(0))+d(\gamma(t),\gamma'(\tau))$.  Since also
$d(\gamma(t),\gamma'(t))\le d(\gamma(t), \gamma'(\tau))+|\tau-t|$ it follows
on substituting for $|\tau-t|$ that
$d(\gamma(t),\gamma'(t))\le2d(\gamma(t),\gamma'(\tau))+d(\gamma(0),\gamma'(0))$.
In particular for divergent rays the distance from a point on one ray to the
other ray tends to infinity with the point.

\begin{corollary} \label{disaxes} A geodesic ray in $\mathcal{T}$ and the axes
of an irreducible mapping class are either asymptotic or divergent.  Two
irreducible mapping classes have axes that coincide or are divergent.
\end{corollary}

\prf  We first consider a ray $\gamma$ and an irreducible element $S$ with
axis $\gamma_S$.  By Proposition \ref{discon}   for unit-speed
parameterizations  the distance between corresponding points of $\gamma$ and
$\gamma_S$ is a convex function.  In particular the distance has a limit
$L_{\infty}=\lim\limits_{t\rightarrow\infty}d(\gamma(t),\gamma'(t))$.  We will
use that the WP geometry along $\gamma_S$ is periodic to show that: either
$L_{\infty}$ is zero or there is a positive lower bound for the convexity of
the distance and thus $L_{\infty}$ is infinite.  

We revisit formula (\ref{secondvar}) for the one-parameter variation $\beta(t,s)$ ($t$ is now the
parameter for the geodesics {\em connecting} $\gamma_S$ to $\gamma$.)  The
integrand of (\ref{secondvar}) is bounded below by the contribution of the
curvature term $-\bigl<R(V,T)T,V\bigr>$, which in turn is non-negative.  The
curve $\beta(t,s')$, $t_0\le t\le t_1$, is a geodesic with initial point on
$\gamma_S$ and length at least $L_{\infty}$.  $T$ is the tangent field of
$\beta(t,s')$ and $V$ is a Jacobi field along $\beta(t,s')$ with initial
vector having unit length.   Now choose $t'$, $t_0<t'\le L_{\infty}$ such that
the neighborhood of radius $t'-t_0$ about a point of $\gamma_S$ is relatively
compact in $\mathcal{T}$.  For the parameter range $t_0\le t\le t'$ the
geodesic segments $\beta(t,s')$, the tangent fields $T$, and the Jacobi fields
$V$ are all modulo the action of $S$ supported on a compact set: the closure
of the neighborhood of a fundamental segment of $\gamma_S$.  Since the WP
curvature is strictly negative on the compact set, we obtain a positive lower
bound for the evaluations, the desired convexity bound for the distance
function.  It follows that $L_{\infty}$ is either zero or infinite.

Consider next irreducible elements $S$, resp. $Q$, with translation lengths
$L_S$, resp. $L_Q$, and axis $\gamma_S$, resp. $\gamma_Q$; assume the axes are
asymptotic in the forward direction.  Choose a reference point $p'$ on
$\gamma_S$.  Given $\epsilon$ positive, choose positive integers $n, m$  such
that $|nL_S-mL_Q|<\epsilon$.  Further choose a positive integer $k_0$ such
that for $k\ge k_0$ the point $p=S^kp'$ on $\gamma_S$ and the corresponding
point  $q$ on $\gamma_Q$ are at distance at most $\epsilon$.  We have then
that $d(S^np,Q^mq)<2\epsilon$; we thus have that $d(Q^{-m}S^np,p)<3\epsilon$.
Since there is a positive lower bound for the distance between distinct points
of the $Mod$ orbit of $p'$, it follows for $\epsilon$ small  that the
transformation $Q^{-m}S^n$ fixes the sequence of points $S^kp'$, $k\ge k_0$.
It follows that $Q^m$ stabilizes $\gamma_S$.  Since the axes are asymptotic,
for $p$ far along $\gamma_S$ the displacement $d(Q^mp,p)$ is close to
$d(Q^mq,q)$; $Q^m$ realizes its translation length on $\gamma_S$; by the
Theorem the axes coincide.    The proof is complete.  

\bibliographystyle{plain}

\end{document}